\documentclass[12pt,oneside,reqno]{amsart}
\usepackage{indentfirst}
\usepackage{graphicx}
\usepackage{mathrsfs}
\usepackage{stmaryrd}
\usepackage{amsfonts}
\usepackage{enumerate,amsmath,amssymb,amsthm}
\usepackage{color}
\usepackage{float}
\newcommand{\dif}{\mathrm{d}}

\newcommand{\be}{\begin{eqnarray}}
\newcommand{\ee}{\end{eqnarray}}
\newcommand{\ce}{\begin{eqnarray*}}
\newcommand{\de}{\end{eqnarray*}}
\newtheorem{theorem}{Theorem}[section]
\newtheorem{lemma}[theorem]{Lemma}
\newtheorem{remark}[theorem]{Remark}
\newtheorem{definition}[theorem]{Definition}
\newtheorem{proposition}[theorem]{Proposition}
\newtheorem{Examples}[theorem]{Examples}
\newtheorem{corollary}[theorem]{Corollary}
\def\e{\varepsilon}

\def\a{\alpha}

\def\d{\delta}

\def\[{{\Big[}}
\def\]{{\Big]}}
\def\<{{\langle}}
\def\>{{\rangle}}
\def\({{\Big(}}
\def\){{\Big)}}

\def\no{\nonumber}
\def\bt{\begin{theorem}}
\def\et{\end{theorem}}
\def\bl{\begin{lemma}}
\def\el{\end{lemma}}
\def\br{\begin{remark}}
\def\er{\end{remark}}
\def\bx{\begin{Examples}}
\def\ex{\end{Examples}}
\def\bd{\begin{definition}}
\def\ed{\end{definition}}
\def\bp{\begin{proposition}}
\def\ep{\end{proposition}}
\def\bc{\begin{corollary}}
\def\ec{\end{corollary}}

\def\cD{{\mathcal D}}

\def\cM{{\mathcal M}}
\def\cN{{\mathcal N}}

\def\cP{{\mathcal P}}

\def\mE{{\mathbb E}}
\def\mF{{\mathbb F}}

\def\mH{{\mathbb H}}

\def\mN{{\mathbb N}}

\def\mP{{\mathbb P}}

\def\mR{{\mathbb R}}
\def\mS{{\mathbb S}}

\def\sA{{\mathscr A}}
\def\sB{{\mathscr B}}
\def\sC{{\mathscr C}}

\def\sF{{\mathscr F}}

\def\sL{{\mathscr L}}

\def\sV{{\mathscr V}}

\def\geq{\geqslant}
\def\leq{\leqslant}
\allowdisplaybreaks[3]

\begin{document}
	
\allowdisplaybreaks
	
\title{Backward multivalued McKean-Vlasov SDEs and associated variational inequalities}
	
\author{Jun Gong and Huijie Qiao$^*$}
	
	
	
\thanks{This work was supported by NSF of China (No.12071071)}
	
\thanks{$*$ Corresponding author: hjqiaogean@seu.edu.cn}
	
	
\date{}
\dedicatory{School of Mathematics,
		Southeast University\\
		Nanjing, Jiangsu 211189, China}

\begin{abstract}
The work concerns a type of backward multivalued  McKean-Vlasov stochastic differential equations. First, we prove the existence and uniqueness of solutions for backward multivalued  McKean-Vlasov stochastic differential equations. Then, it is presented that their solutions depend continuously on the terminal values. Finally, we give a probabilistic interpretation for viscosity solutions of nonlocal quasi-linear parabolic variational inequalities.
\end{abstract}

\maketitle \rm

\section{Introduction}
McKean-Vlasov stochastic differential equations (SDEs for short), also called mean-field SDEs or distribution-dependent SDEs, can track back to Kac \cite{k1} in 1956. Then in 2009 Buckdahn, Djehiche, Li and Peng \cite{bdlp} investigated a type of backward McKean-Vlasov SDEs. Since this, the theory of backward McKean-Vlasov SDEs and forward-backward McKean-Vlasov SDEs, as well as that of the associated partial differential equations have been widely studied. For example, in \cite{BUCK} Buckdahn, Li and Peng  not only proved the existence, uniqueness and a comparison theorem of the solutions for backward McKean-Vlasov SDEs but also obtained a probabilistic interpretation of related nonlocal partial differential equations. Later, Li and Luo \cite{ll} studied reflected backward McKean-Vlasov SDEs, and proved the existence and uniqueness for their solutions under the Lipschitz condition. Li \cite{L1} also observed reflected backward McKean-Vlasov SDEs in a purely probabilistic method, and gave a probabilistic interpretation for obstacle problems of nonlinear and nonlocal partial differential equations by means of reflected backward McKean-Vlasov SDEs. Besides, Lu, Ren and Hu \cite{lrh} dealt with a class of backward McKean-Vlasov SDEs with subdifferential operators corresponding to lower semi-continuous convex functions. By means of the Yosida approximation, they established the existence and uniqueness of the solutions for backward McKean-Vlasov SDEs with subdifferential operators, and attained a probability interpretation for the viscosity solutions of a class of nonlocal parabolic variational inequalities. Here we emphasize that in \cite{BUCK, L1, ll, lrh} the coefficients of backward McKean-Vlasov SDEs depend on the distributions of solution processes through their expectations.

On the other hand, there also exist backward McKean-Vlasov SDEs in which the coefficients straightly depend on the distributions of solution processes. Note that the type of backward McKean-Vlasov SDEs is more general than the mentioned type in the above paragraph. Let us list some works related with ours. In \cite{cd} Carmona and Delarue provided an existence result for the solutions of a fully coupled forward-backward McKean-Vlasov SDE under a very mild Lipschitz condition. Then in \cite{byz} Bensoussan, Yam and Zhang proposed a broad class of natural monotonicity conditions and established the well-posedness for a type of forward-backward McKean-Vlasov SDEs. Recently, Li \cite{L} studied a type of backward McKean-Vlasov SDEs driven by Brownian motions and independent Poisson random measures. There she proved the well-posedness of their solutions, and also obtained the well-posedness of classical solutions for related nonlocal quasi-linear integral-partial differential equations under some regular assumptions.

In the paper, we concentrate on a type of more general backward McKean-Vlasov SDEs. Concretely speaking, we fix $T>0$ and a complete filtered probability space $(\Omega,\mathscr{F},\{\mathscr{F}_t\}_{t\in[0,T]},\mP)$, and consider the following backward multivalued McKean-Vlasov SDE on $\mR^d$:
\be\left\{\begin{array}{l}
\dif Y_t \in \ A(Y_t)\dif t -G(t,Y_t,Z_t,\sL_{(Y_t, Z_t)})\dif t + Z_t\dif W_t,\\
Y_T=\xi, \sL_{(Y_t, Z_t)}=$the probability  distribution of~$(Y_t, Z_t),
\end{array}
\label{eq1}
\right.
\ee
where $W_{\cdot}=(W_{\cdot}^1,W_{\cdot}^2,\cdots,W_{\cdot}^l)$ is a $(\sF_t)_{t\in[0,T]}$-adapted $l$-dimensional standard Brownian motion, $A:\mR^d \mapsto 2^{\mR^d} $ is a maximal monotone operator, for the coefficient $G: \Omega\times [0,T] \times \mR^d\times \mR^{d\times l}\times\cM_2(\mR^d\times\mR^{d\times l})\mapsto{\mR^d}$, $\forall (y,z,\vartheta)\in\mR^d\times \mR^{d\times l}\times\cM_2(\mR^d\times\mR^{d\times l})$, $G(\cdot,y,z,\vartheta)$ is $(\sF_t)_{t\in[0,T]}$-predictable, and $\xi$ is a $\sF_T$-measurable random variable with values in $\overline{\cD\left( A \right) }$ (See Subsection \ref{mmo}) and $\mE|\xi|^{2}<\infty$. If $A=0$, Eq.(\ref{eq1}) becomes a backward McKean-Vlasov SDE in \cite{byz, cd, L}. And if $A\neq 0$, backward multivalued McKean-Vlasov SDEs like Eq.(\ref{eq1}) include reflected backward McKean-Vlasov SDEs in \cite{L1, ll} and backward McKean-Vlasov SDEs with subdifferential operators in \cite{lrh}. Hence, as far as we know, backward multivalued McKean-Vlasov SDEs are the most general backward McKean-Vlasov SDEs. Since backward multivalued McKean-Vlasov SDEs like Eq.(\ref{eq1}) are widely applied in finance, the control theory and the game theory (cf. \cite{drh} and cited references there), we are devoted to studying Eq.(\ref{eq1}). 

As a whole, our contribution are two-folded:

$\bullet$ We prove the well-posedness of Eq.(\ref{eq1}) under the Lipschitz condition, and also obtain that its solution depend continuously on the terminal value.

$\bullet$ Combining Eq.(\ref{eq1}) with a forward McKean-Vlasov SDE, we give a probabilistic interpretation for viscosity solutions of nonlocal quasi-linear parabolic variational inequalities. 

It is worthwhile to mentioning our results and methods. Notice that Li \cite{L1} and Lu, Ren, Hu \cite{lrh} also showed the well-posedness and probabilistic interpretation. Since our equation is more general than that in \cite{L1, lrh}, our results are better. Moreover, we establish continuous dependence of solutions on the terminal values. This is important for application of backward multivalued McKean-Vlasov SDEs. Besides, we emphasize that the appearance of maximal monotone operators brings a lot of trouble, and our methods and techniques are more subtle.

The content of the paper is arranged as follows. In the next section, we introduce notations and concepts, such as maximal monotone operators and the derivative for functions on $\cM_{2}(\mR^d)$. Moreover, a result about backward McKean-Vlasov SDEs is listed and some conclusions on backward multivalued stochastic differential equations are proved in the section. Then we place the existence, uniqueness and continuous dependence on the terminal values of the solutions for Eq.(\ref{eq1}) in Section \ref{teus}. In Section \ref{CONN}, viscosity solutions of nonlocal quasi-linear parabolic variational inequalities are established.

The following convention will be used throughout the paper: $C$ with or without indices will denote different positive constants whose values may change from one place to another.

\section{Preliminary}\label{fram}

In the section, we introduce notations and concepts, and recall some results used in the sequel.

\subsection{Notations}\label{nota}

In the subsection, we introduce some notations.

For convenience, we shall use $\mid\cdot\mid$ and $\parallel\cdot\parallel$  for norms of vectors and matrices, respectively. Furthermore, let $\langle\cdot$ , $\cdot\rangle$ denote the scalar product in $\mR^d$. Let $B^*$ denote the transpose of a matrix $B$.

Let $\sB(\mR^d)$ be the Borel $\sigma$-algebra on $\mR^d$ and $\cP({\mR^d})$ be the space of all probability measures defined on $\sB(\mR^d)$ carrying the usual topology of weak convergence. Let $\cM_{2}(\mR^d)$ be the set of probability measures on $\sB(\mR^d)$ with finite second order moments. That is,
$$
\cM_2\left( \mathbb{R}^d \right) :=\left\{ \mu \in \cP\left( \mathbb{R}^d \right) :\left\| \mu \right\| _{2}^{2}:=\int_{\mathbb{R}^d}{\left| x \right|^2\mu \left( \dif x \right) <\infty} \right\}. 
$$
Define  the following metric on $\cM_{2}(\mR^d)$:
\ce
\rho^2(\mu_1,\mu_2):=\inf_{\pi\in\sC(\mu_1, \mu_2)}\int_{\mR^d\times\mR^d}|x-y|^2\pi(\dif x, \dif y), \quad \mu_1, \mu_2\in\cM_2(\mR^d),
 \de
where $\sC(\mu_1, \mu_2)$ denotes the set of  all the probability measures whose marginal distributions are $\mu_1, \mu_2$, respectively. Thus, $(\cM_2(\mR^d), \rho)$ is a Polish space.

Let $\mS^2_{\mF}([0,T],\mR^d)$ be the space of all $\mF$-adapted processes $Y: \Omega\times [0,T]\mapsto\mR^d$ with $\mE\left[\sup\limits_{t\in[0,T]}|Y_t|^2\right]<\infty$, where $\mF:=(\sF_t)_{t\in[0,T]}$. Let $\mH_{\mF}^2([0,T],\mR^{d\times l})$ be the space of all $\mF$-predictable processes $Z: \Omega\times [0,T]\mapsto\mR^{d\times l}$ with $\mE\left[\int_0^T\|Z_t\|^2\dif t\right]<\infty$.

Let $C(\mR^d)$ be the collection of continuous functions on $\mR^d$ and $C^2(\mR^d)$ be the space of continuous functions on $\mR^d$ which have continuous partial derivatives of order up to $2$.

\subsection{Maximal monotone operators}\label{mmo}

In the subsection, we introduce maximal monotone operators. 

For a multivalued operator $A: \mR^d\mapsto 2^{\mR^d}$, where $2^{\mR^d}$ stands for all the subsets of $\mR^d$, set
\ce
&&\cD(A):= \left\{x\in \mR^d: A(x) \ne \emptyset\right\},\\
&&Gr(A):= \left\{(x,y)\in \mR^{2d}:x \in \cD(A), ~ y\in A(x)\right\}.
\de
We say that $A$ is monotone if $\langle x_1 - x_2, y_1 - y_2 \rangle \geq 0$ for any $(x_1,y_1), (x_2,y_2) \in Gr(A)$, and $A$ is maximal monotone if 
$$
(x_1,y_1) \in Gr(A) \iff \langle x_1-x_2, y_1 -y_2 \rangle \geq 0, \forall (x_2,y_2) \in Gr(A).
$$

In the following, we recall some properties of a maximal monotone operator $A$. (cf.\cite{cepa1, gp})

\begin{enumerate}[(i)]
\item
${\rm Int}(\cD(A))$ and $\overline{\mathrm{\cD}(A)}$ are convex subsets of $\mR^d$ with ${\rm Int}\left( \overline{\mathrm{\cD}(A)} \right) = {\rm Int}\( \mathrm{\cD}(A) \) 
$, where ${\rm Int}(\cD(A))$ denotes the interior of the set $\cD(A)$. 
\item For every $x\in\mR^d$, $A(x)$ is a closed and convex subset of $\mR^d$. Let $A^{\circ}\left( x \right):= proj_{A(x)}(0)$ be the minimal section of $A$, where $proj_D$ is designated as the projection on every closed and
convex subset $D\subset\mR^d$ and $proj_{\emptyset}(0) =\infty$. Then
$$
x\in\cD(A) \Longleftrightarrow |A^{\circ}\left( x \right)|<\infty.
$$
\item For $\varepsilon > 0$, the resolvent operator $J_{\varepsilon}:=\left( I+\varepsilon A \right) ^{-1}$ is a single-valued and contractive operator defined on $\mR^d$ and takes values in $\mathrm{\cD}(A)$, and 
$$
\lim\limits_{\e\downarrow 0}J_\e(x)=proj_{\overline{\mathrm{\cD}(A)}}(x), \quad x\in\mR^d.
$$
\item $A_{\varepsilon}:= \frac{1}{\varepsilon}\left(I-J_{\varepsilon}\right)$, called the Yosida approximation of $A$, is also a single-valued, maximal monotone and Lipschitz continuous operator with the Lipschitz constant $\frac{1}{\varepsilon}$. 
\item 
$A_{\varepsilon}\left( x \right) \,\,\in \,\,A\left( J_{\varepsilon}\left( x \right) \right), \quad x\in\mR^d$.
\item 
$|A_{\varepsilon}\left( x \right) | \leq |A^{\circ}\left( x \right)|, \quad x\in\cD(A)$.
\item
For any $x\in\cD(A)$, $\lim\limits_{\varepsilon \downarrow 0}A_{\varepsilon}\left( x \right) =A^{\circ}\left( x \right)$, and
\ce
&&\lim_{\varepsilon \downarrow 0}|A_{\varepsilon}\left( x \right) |=|A^{\circ}\left( x \right)|, \quad x\in\cD(A),\\
&&\lim_{\varepsilon \downarrow 0}|A_{\varepsilon}\left( x \right) |=\infty,\qquad\qquad x\notin\cD(A).
\de
\end{enumerate}

About the Yosida approximation $A_{\varepsilon}$, we also mention the following property (\cite[Lemma 5.4]{cepa1}). 

\bl\label{yosi}
There exist three constants $a\in\mR^d, M_1>0, M_2\geq 0$ only dependent on $A$ such that for any $\e>0$ and $x\in\mR^d$
$$
\<A_{\varepsilon}(x), x-a\>\geq M_1|A_{\varepsilon}(x)|-M_2|x-a|-M_1M_2.
$$
\el

\medspace

Let $\sV_{0}$ be the set of all continuous functions $K: [0,T]\mapsto\mR^d$ with finite variations and $K_{0} = 0$. For $K\in\sV_0$ and $s\in [0,T]$, we shall use $|K|_{0}^{s}$ to denote the variation of $K$ on [0,s]
and write $|K|_{TV}:=|K|_{0}^{T}$. Set
\ce
&&\sA:=\Big\{(Y,K): Y\in C([0,T],\overline{\cD(A)}), K \in \sV_0, \\
&&\qquad\qquad\quad~\mbox{and}~\langle Y_{t}-x, \dif K_{t}-y\dif t\rangle \geq 0 ~\mbox{for any}~ (x,y)\in Gr(A)\Big\}.
\de
And about $\sA$ we have two following results (cf.\cite{cepa2, ZXCH}).

\bl\label{equi}
For $Y\in C([0,T],\overline{\cD(A)})$ and $K\in \sV_{0}$, the following statements are equivalent:
\begin{enumerate}[(i)]
	\item $(Y,K)\in \sA$.
	\item For any $(x,y)\in C([0,T],\mR^d)$ with $(x_t,y_t)\in Gr(A)$, it holds that 
	$$
	\left\langle Y_t-x_t, \dif K_t-y_t\dif t\right\rangle \geq0.
	$$
	\item For any $(Y^{'},K^{'})\in \sA$, it holds that 
	$$
	\left\langle Y_t-Y_t^{'},\dif K_t-\dif K_t^{'}\right\rangle \geq0.
	$$
\end{enumerate}
\el

\bl\label{limiconv}
Assume that $\{K^{n},n\in \mN\}\subset \sV_{0}$ converges to some $K$ in $C([0,T];\mR^d)$ and $\underset{n\in \mathbb{N}}{\sup}\left| K^n \right|_{TV}<\infty $. Then $K\in \sV_0$, and 
$$
\underset{n\rightarrow \infty}{\lim}\int_0^T{\left< Y_{s}^{n}, \dif K_{s}^{n} \right>}=\int_0^T{\left< Y_{s}, \dif K_{s}\right>},
$$
where the sequence $\{Y^{n}\}\subset C([0,T];\mR^d)$ converges to some $Y$ in $C([0,T];\mR^d)$.
\el

\subsection{Backward McKean-Vlasov SDEs}\label{bmvsde}

In the subsection, we recall a result about backward McKean-Vlasov SDEs. 

Now given a function $\tilde{G}: \Omega\times [0,T] \times \mR^d\times \mR^{d\times l}\times\cM_2(\mR^d\times\mR^{d\times l})\mapsto{\mR^d}$ satisfying that for $\forall (y,z,\vartheta)\in\mR^d\times \mR^{d\times l}\times\cM_2(\mR^d\times\mR^{d\times l})$, $\tilde{G}(\cdot,y,z,\vartheta)$ is $(\sF_t)_{t\in[0,T]}$-predictable, and a $\sF_T$-measurable random variable $\xi$ with $\mE|\xi|^{2}<\infty$. 

Assume:
\begin{enumerate}[(${\bf H}^1_{\tilde{G}}$)]
\item There exsits a non-random constant $L_{\tilde{G}}>0$ such that for a.s. $\omega\in\Omega$, it holds that
\ce
&&\left| \tilde{G}\left( \omega ,t,y,z,\vartheta  \right) -\tilde{G}\left( \omega ,t,y^{\prime},z^{\prime},\vartheta ^{\prime} \right) \right|\leq L_{\tilde{G}}\left( \left| y-y^{\prime} \right|+\| z-z^{\prime} \|+\rho \left( \vartheta ,\vartheta ^{\prime} \right) \right),\\
 &&t\in\left[ 0, T \right], y, y^{'}\in \mR^d, z, z^{'}\in \mR^{d\times l}, \vartheta, \vartheta^{'}, \in \cM_{2}(\mR^d\times\mR^{d\times l}).
\de
\end{enumerate}
\begin{enumerate}[(${\bf H}^2_{\tilde{G}}$)]
\item $\mathbb{E}\int_0^T{\left| \tilde{G}\left( t,0,0,\delta _{({\bf 0}, {\bf 0})} \right) \right|}^2\dif t<\infty $, where $\delta_{({\bf 0}, {\bf 0})} $ is the Dirac measure at $({\bf 0}, {\bf 0})$.
\end{enumerate}

The following result comes from \cite[Theorem A.1.]{L}.

\bt\label{bmvsdeeu}
Assume that (${\bf H}^1_{\tilde{G}}$)-(${\bf H}^2_{\tilde{G}}$) hold. Then the following backward McKean-Vlasov SDE
\ce
\left\{ \begin{array}{l}
	\mathrm{d}Y_{t}=-\tilde{G}\(t, Y_{t}, Z_{t}, \sL_{(Y_{t}, Z_{t})}\)\dif t+Z_{t}\dif W_t, \quad 0\leq t<T,\\
	Y_{T}=\xi,
\end{array}
\right. 
\de 
has a unique solution $(Y_{\cdot}, Z_{\cdot})\in \mS_{\mF}^2([0,T],\mR^d)\times\mH_{\mF}^2([0,T],\mR^{d\times l})$. 
\et

\subsection{Backward multivalued SDEs}\label{bmsde}

In the subsection, we recall and prove some results about backward multivalued SDEs. 

Consider the following backward multivalued SDE on $\mR^d$:
\be\left\{\begin{array}{l}
\dif Y_t \in \ A(Y_t)\dif t -g(t,Y_t,Z_t)\dif t + Z_t\dif W_t,\\
Y_T=\xi, 
\end{array}
\right.
\label{bmsdeeq}
\ee
where the coefficient $g: \Omega\times [0,T] \times \mR^d\times \mR^{d\times l}\mapsto{\mR^d}$ is Borel measurable, $\forall (y,z)\in\mR^d\times \mR^{d\times l}$, $g(\cdot,y,z)$ is $(\sF_t)_{t\in[0,T]}$-progressively measurable, and $\xi$ is a $\sF_T$-measurable random variable with values in $\overline{\cD\left( A \right) }$ and $\mE|\xi|^{2}<\infty$. We define the existence and uniqueness of solutions for Eq.(\ref{bmsdeeq}).

\bd\label{solu}
We say that Eq.$(\ref{bmsdeeq})$ admits a solution with the terminal value $\xi$ if there exists a triple $\left\{ \left( Y_t, K_t, Z_t\right) :t\in \left[ 0,T \right] \right\} $ which is a $(\sF_t)_{t\in[0,T]}$-progressively measurable process and satisfies 
\begin{enumerate}[(i)]
\item 
$(Y_{\cdot},K_{\cdot})\in \sA$, $\dif \mP\times\dif t$-a.e. on $\Omega\times\left[ 0, T \right]$,
\item 
$$
\mathbb{E}\left( \underset{0\leqslant t\leqslant T}{\sup}\left| Y_t \right|^2+\int_0^T{\| Z_s \|^2}\dif s \right) <\infty, 
$$
\item 
\ce
Y_t=\xi -\left( K_T-K_t \right)+\int_t^T{g\left( s, Y_s, Z_s\right)}\dif s-\int_t^T{Z_s}\dif W_s.
 \de
\end{enumerate}
\ed

\bd\label{uniq}
Suppose that $\{(\Omega,\mathscr{F},\{\mathscr{F}_t\}_{t\in[0,T]},\mP), (Y^1, K^1, Z^1)\}$ and $\{(\Omega,\mathscr{F},\{\mathscr{F}_t\}_{t\in[0,T]},\mP),\\ (Y^2, K^2, Z^2)\}$ are two solutions for Eq.(\ref{bmsdeeq}) with $Y^1_{T}=Y^2_T=\xi, K^1_{T}=K^2_T$. If $Y^1_t=Y^2_t, 
K^1_t=K^2_t, t\in[0,T], \mP$-a.s. and $Z^1_{\cdot}=Z^2_{\cdot}, \dif \mP\times\dif t$-a.e., we say that the uniqueness holds for Eq.(\ref{bmsdeeq}).
\ed

For the existence and uniqueness of Eq.(\ref{bmsdeeq}), we assume:
\begin{enumerate}[(${\bf H}^1_{g}$)]
\item There exsits a non-random constant $L_{g}>0$ such that for a.s. $\omega\in\Omega$, it holds that
\ce
&&\left| g\left( \omega ,t,y,z\right) -g\left( \omega ,t,y^{\prime},z^{\prime}\right) \right|\leq L_{g}\left( \left| y-y^{\prime} \right|+\| z-z^{\prime} \|\right),\\
 &&t\in\left[ 0, T \right], y, y^{'}\in \mR^d, z, z^{'}\in \mR^{d\times l}.
\de
\end{enumerate}
\begin{enumerate}[(${\bf H}^2_{g}$)]
\item $\mathbb{E}\int_0^T{\left| g\left( t,0,0\right) \right|}^2\dif t<\infty $.
\end{enumerate}
\begin{enumerate}[(${\bf H}_{A}$)]
\item $\cD(A)=\mR^d$ and $|A^\circ(x)|\leq C(1+|x|), x\in\mR^d$, where $C>0$ is a constant.
\end{enumerate}

\br
If $A=\partial \varphi$, where $\varphi$ is a proper ($\varphi\not\equiv+\infty$), convex and lower semicontinuous function, $\varphi(x)\geq \varphi(0)=0$ and $\partial \varphi$ is the subdifferential operator of $\varphi$, $({\bf H}_{A})$ can be replaced by $\mE|\varphi(\xi)|<\infty$.
\er

The following result comes from \cite[Theorem 2 and Theorem 7]{MODE}.

\bt\label{bmsdeeu}
Assume that $({\bf H}_{A})$ holds and the coefficient $g$ satisfies $({\bf H}^1_g)$-$({\bf H}^2_g)$. Then there exists a unique solution for Eq.(\ref{bmsdeeq}).
\et

Next, suppose that $(Y^i,K^i,Z^i)$ is the unique solution of the following backward multivalued SDE:
\ce\left\{\begin{array}{l}
\dif Y^i_t \in \ A(Y^i_t)\dif t -g^i(t,Y^i_t,Z^i_t)\dif t + Z^i_t\dif W_t,\\
Y^i_T=\xi^i, 
\end{array}
\right.
\de
where the coefficient $g^i: \Omega\times [0,T] \times \mR^d\times \mR^{d\times l}\mapsto{\mR^d}$ is Borel measurable, $\forall (y,z)\in\mR^d\times \mR^{d\times l}$, $g^i(\cdot,y,z)$ is $(\sF_t)_{t\in[0,T]}$-progressively measurable, and $\xi^i$ is a $\sF_T$-measurable random variable with values in $\overline{\cD\left( A \right) }$ and $\mE|\xi^i|^{2}<\infty$, $i=1,2$. We have the following result.

\bt\label{diffesti}
Assume that $({\bf H}_{A})$ holds, $g^1, g^2$ satisfy $({\bf H}^1_g)$-$({\bf H}^2_g)$, and there exists a constant $c^*>0$ such that for a.s. $\omega\in\Omega$, it holds that
\ce
&&\left| g^i\left( \omega ,t,y,z\right) -g^i\left( \omega ,t,y^{\prime},z^{\prime}\right) \right|\leq c^*\left( \left| y-y^{\prime} \right|+\| z-z^{\prime} \|\right),\\
 &&t\in\left[ 0, T \right], y, y^{'}\in \mR^d, z, z^{'}\in \mR^{d\times l}, i=1,2.
\de
Then for $\tilde{Y}:=Y^1-Y^2, \tilde{Z}:=Z^1-Z^2$, $\tilde{g}:=g^1-g^2, \tilde{\xi}:=\xi^1-\xi^2$, it holds that

(i) For $\d>0$, there exists a constant $\a>0$ such that for any $t\in[0,T]$
\be
&&|\tilde{Y}_t|^2+\frac{1}{2}\mE\left[\int_t^Te^{\a(s-t)}(|\tilde{Y}_s|^2+\|\tilde{Z}_s\|^2)\dif s|\sF_t\right]\no\\
&\leq&\mE[e^{\a(T-t)}|\tilde{\xi}|^2|\sF_t]+c^*\d\mE\left[\int_t^Te^{\a(s-t)}|\tilde{g}(s,Y_s^1,Z_s^1)|^2\dif s|\sF_t\right], a.s..
\label{diffesti1}
\ee

(ii) There exists a constant $C>0$ such that for any $t\in[0,T]$
\be
\mE\left[\sup\limits_{s\in[t,T]}|\tilde{Y}_s|^2+\int_t^T\|\tilde{Z}_s\|^2\dif s|\sF_t\right]\leq C\mE\left[|\tilde{\xi}|^2+\int_t^T|\tilde{g}(s,Y_s^1,Z_s^1)|^2\dif s|\sF_t\right], a.s..
\ee
\et
\begin{proof}
First of all, we prove (i). Applying It\^{o}'s formula to $e^{\alpha t}|\tilde{Y}_t|^2$, where $\a$ is determined later, by Lemma \ref{equi}, we get that
\ce
&&e^{\alpha t}\left| \tilde{Y} _t\right|^2+\int_t^T{e^{\alpha s}\left\| \tilde{Z}_s \right\|}^2\dif s+\alpha \int_t^T{e^{\alpha s}\left| \tilde{Y_s} \right|^2}\dif s
\\
&=&e^{\alpha T}\left| \tilde{\xi } \right|^2+2\int_t^T{e^{\alpha s}\left< \tilde{Y_s},g^1\left( s,Y_{s}^{1},Z_{s}^{1} \right) -g^2\left( s,Y_{s}^{2},Z_{s}^{2} \right) \right>}\dif s\\
&&-2\int_t^T{e^{\alpha s}\left< \tilde{Y}_s,\tilde{Z_s}\dif W_s \right>}-2\int_t^T{e^{\alpha s}}\left< \tilde{Y}_s,\dif \tilde{K_s} \right> \\
&\leq&e^{\alpha T}\left| \tilde{\xi } \right|^2+2\int_t^T{e^{\alpha s}\left< \tilde{Y}_s,\tilde{g}\left( s,Y_{s}^{1},Z_{s}^{1} \right) +g^2\left( s,Y_{s}^{1},Z_{s}^{1} \right) -g^2\left( s,Y_{s}^{2},Z_{s}^{2} \right) \right>}\dif s\\
&&-2\int_t^T{e^{\alpha s}\left< \tilde{Y}_s,\tilde{Z_s}\dif W_s \right>}.
\de
Taking the conditional expectation with respect to $\sF_t$ on two sides of the above inequality, we have that
\ce
&&e^{\alpha t}\left| \tilde{Y}_t \right|^2+\mathbb{E}\left[ \int_t^T{e^{\alpha s}\left\| \tilde{Z}_s \right\|}^2ds+\alpha \int_t^T{e^{\alpha s}\left| \tilde{Y}_s \right|^2}\dif s|\sF_t\right]
\\
&\leqslant& \mathbb{E}\left[ e^{\alpha T}\left| \tilde{\xi } \right|^2|\sF_t \right] +2\mathbb{E}\left[ \int_t^T{e^{\alpha s}\left< \tilde{Y_s},\tilde{g}\left( s,Y_{s}^{1},Z_{s}^{1} \right) +g^2\left( s,Y_{s}^{1},Z_{s}^{1} \right) -g^2\left( s,Y_{s}^{2},Z_{s}^{2} \right) \right>}\dif s|\sF_t \right].
\de
By $({\bf H}_{g}^{1})$ and the Young inequality, it holds that
\ce
&&e^{\alpha t}\left| \tilde{Y}_t \right|^2+\mathbb{E}\left[ \int_t^T{e^{\alpha s}\left\| \tilde{Z}_s \right\|}^2ds+\alpha \int_t^T{e^{\alpha s}\left| \tilde{Y}_s \right|^2}\dif s|\sF_t\right]\\
&\leqslant& \mathbb{E}\left[ e^{\alpha T}\left| \tilde{\xi } \right|^2|\sF_t \right] +2\mathbb{E}\left[ \int_t^T e^{\alpha s}|\tilde{Y_s}|\cdot |\tilde{g}\left( s,Y_{s}^{1},Z_{s}^{1} \right)|\dif s|\sF_t \right]\\
&&+2\mathbb{E}\left[ \int_t^T e^{\alpha s}|\tilde{Y_s}|\cdot |g^2\left( s,Y_{s}^{1},Z_{s}^{1} \right) -g^2\left( s,Y_{s}^{2},Z_{s}^{2} \right)|\dif s|\sF_t \right]\\
&\leqslant& \mathbb{E}\left[ e^{\alpha T}\left| \tilde{\xi } \right|^2|\sF_t \right] +c^*\delta \mathbb{E}\left[\int_t^T{e^{\alpha s}\left| \tilde{g}\left( s,Y_{s}^{1},Z_{s}^{1} \right) \right|}^2ds|\sF_t\right]+\frac{1}{c^*\delta}\mathbb{E}\left[\int_t^T{e^{\alpha s}\left| \tilde{Y_s} \right|^2}ds|\sF_t\right]\\
&&+2c^*\mathbb{E}\left[\int_t^T{e^{\alpha s}\left| \tilde{Y_s} \right|^2}\dif s|\sF_t\right]+\frac{1}{2}\mathbb{E}\left[\int_t^T{e^{\alpha s}\left\| \tilde{Z_s} \right\| ^2}\dif s|\sF_t\right]+2c^{*2}\mathbb{E}\left[\int_t^T{e^{\alpha s}\left| \tilde{Y_s} \right|^2}\dif s|\sF_t\right]\\
&\leqslant& \mathbb{E}\left[ e^{\alpha T}\left| \tilde{\xi } \right|^2|\sF_t  \right] +c^*\delta \mathbb{E}\left[\int_t^T{e^{\alpha s}\left| \tilde{g}\left( s,Y_{s}^{1},Z_{s}^{1} \right) \right|}^2ds|\sF_t\right]+\frac{1}{2}\mathbb{E}\left[\int_t^T{e^{\alpha s}\left\| \tilde{Z_s} \right\| ^2}\dif s|\sF_t\right]\\
&&+\left( \frac{1}{c^*\delta}+2c^*+2c^{*2} \right) \mathbb{E}\left[\int_t^T{e^{\alpha s}\left| \tilde{Y_s} \right|^2}\dif s|\sF_t\right],
\de
where $\d>0$ is a constant. Taking $\alpha >\left( \frac{1}{c^*\delta}+2c^*+2c^{*2}+\frac{1}{2} \right)$, we obtain (\ref{diffesti1}).

Next, we prove (ii). For any $t\in[0,T]$, making use of the It\^{o} formula and Lemma \ref{equi}, we attain that for $s\in[t,T]$    
\be
&&\left| \tilde{Y}_s \right|^2+\int_s^T{\left\| \tilde{Z} _r\right\|}^2\dif r\no\\
&=&\left| \tilde{\xi } \right|^2+2\int_s^T{\left< \tilde{Y}_r,\tilde{g}\left( r,Y_{r}^{1},Z_{r}^{1} \right) +g^2\left( r,Y_{r}^{1},Z_{r}^{1} \right) -g^2\left( r,Y_{r}^{2},Z_{r}^{2} \right) \right>}\dif r\no\\
&&-2\int_s^T{\left< \tilde{Y}_r,\tilde{Z}_r\dif W_r \right>}-2\int_s^T{\left< \tilde{Y}_r, \dif \tilde{K}_r \right>}\no\\
&\leq&\left| \tilde{\xi } \right|^2+2\int_s^T{\left< \tilde{Y}_r,\tilde{g}\left( r,Y_{r}^{1},Z_{r}^{1} \right) +g^2\left( r,Y_{r}^{1},Z_{r}^{1} \right) -g^2\left( r,Y_{r}^{2},Z_{r}^{2} \right) \right>}\dif r\no\\
&&-2\int_s^T{\left< \tilde{Y}_r,\tilde{Z}_r\dif W_r \right>}.
\label{tilest}
\ee    
Taking the conditional expectation with respect to $\sF_t$, by $({\bf H}_{g}^{1})$, we get 
\ce
&&\mathbb{E}\left[ \left| \tilde{Y}_s \right|^2|\sF_t\right]+\mE\left[\int_s^T{\left\| \tilde{Z}_r \right\|}^2\dif r|\sF_t \right] \no
\\
&\leqslant& \mE\left[\left| \tilde{\xi } \right|^2|\sF_t\right]+2\mE\left[\int_s^T\left|\left< \tilde{Y}_r,\tilde{g}\left( r,Y_{r}^{1},Z_{r}^{1} \right) \right>\right|\dif r|\sF_t \right]\no\\
&&+2\mE\left[\int_s^T\left|\left< \tilde{Y}_r,g^2\left( r,Y_{r}^{1},Z_{r}^{1} \right) -g^2\left( r,Y_{r}^{2},Z_{r}^{2} \right) \right>\right|\dif r|\sF_t \right]\no\\
&\leqslant& \mE\left[\left| \tilde{\xi } \right|^2|\sF_t\right]+\mE\left[\int_s^T|\tilde{Y}_r|^2\dif r|\sF_t \right]+\mE\left[\int_s^T|\tilde{g}\left( r,Y_{r}^{1},Z_{r}^{1} \right) |^2\dif r|\sF_t \right]\no\\
&&+2c^*\mE\left[\int_s^T|\tilde{Y}_r|^2\dif r|\sF_t \right]+2c^{*2}\mE\left[\int_s^T|\tilde{Y}_r|^2\dif r|\sF_t \right]+\frac{1}{2}\mE\left[\int_s^T{\left\| \tilde{Z}_r \right\|}^2\dif r|\sF_t \right]\no\\
&\leqslant& \mE\left[\left| \tilde{\xi } \right|^2|\sF_t\right]+\mE\left[\int_s^T|\tilde{g}\left( r,Y_{r}^{1},Z_{r}^{1} \right) |^2\dif r|\sF_t \right]+(1+2c^*+2c^{*2})\int_s^T\mE\left[|\tilde{Y}_r|^2|\sF_t \right]\dif r\no\\
&&+\frac{1}{2}\mE\left[\int_s^T{\left\| \tilde{Z}_r \right\|}^2\dif r|\sF_t \right],
\de
which together with the Gronwall inequality yields that
\be
\sup\limits_{s\in[t,T]}\mathbb{E}\left[ \left| \tilde{Y}_s \right|^2|\sF_t\right]+\mE\left[\int_t^T{\left\| \tilde{Z}_r \right\|}^2\dif r|\sF_t \right] \leq C\mE\left[\left| \tilde{\xi } \right|^2|\sF_t\right]+\mE\left[\int_t^T|\tilde{g}\left( r,Y_{r}^{1},Z_{r}^{1} \right) |^2\dif r|\sF_t \right].
\label{iequ1}
\ee
	    
Now, we again observe (\ref{tilest}). By the Young inequality, the BDG inequality and $({\bf H}_{g}^{1})$, it holds that
\ce
&&\mathbb{E}\left[ \underset{s\in \left[ t,T \right]}{\sup}\left| \tilde{Y_s} \right|^2+\int_t^T{\left\| \tilde{Z}_r \right\|}^2\dif r\left| \mathscr{F}_t \right. \right] \\
&\leq&\mathbb{E}\left[ \left| \tilde{\xi } \right|^2\left| \mathscr{F}_t \right. \right] +2\mE\left[\int_t^T{\left|\left< \tilde{Y}_r,\tilde{g}\left( r,Y_{r}^{1},Z_{r}^{1} \right) +g^2\left( r,Y_{r}^{1},Z_{r}^{1} \right) -g^2\left( r,Y_{r}^{2},Z_{r}^{2} \right) \right>\right|}\dif r|\sF_t\right] \no\\
&&+2\mE\left[\underset{s\in \left[ t,T \right]}{\sup}\left|\int_s^T{\left< \tilde{Y}_r,\tilde{Z}_r\dif W_r \right>}\right||\sF_t\right]\\
&\leqslant& \mathbb{E}\left[ \left| \tilde{\xi } \right|^2\left| \mathscr{F}_t \right. \right] +\mE\left[\int_t^T|\tilde{Y}_r|^2\dif r|\sF_t \right]+\mE\left[\int_t^T|\tilde{g}\left( r,Y_{r}^{1},Z_{r}^{1} \right) |^2\dif r|\sF_t \right] \\
&&+2c^*\mathbb{E}\left[ \int_t^T\left| \tilde{Y}_r \right|^2\dif r\left| \mathscr{F}_t \right. \right] +4c^{*2}\mathbb{E}\left[ \int_t^T\left| \tilde{Y}_r\right|^2\dif r\left| \mathscr{F}_t \right. \right] \\
&&+\frac{1}{4}\mathbb{E}\left[ \int_t^T{\left\| \tilde{Z} _r\right\| ^2}\dif r\left| \mathscr{F}_t \right. \right] +C\mathbb{E}\left[ \left( \int_t^T{\left| \tilde{Y}_r \right|^2\left\| \tilde{Z}_r\right\| ^2}\dif r \right) ^{\frac{1}{2}}\left| \mathscr{F}_t \right. \right] \\
&\leqslant& \mathbb{E}\left[ \left| \tilde{\xi } \right|^2\left| \mathscr{F}_t \right. \right]+\mE\left[\int_t^T|\tilde{g}\left( r,Y_{r}^{1},Z_{r}^{1} \right) |^2\dif r|\sF_t \right]  +(1+2c^*+4c^{*2})\int_t^T\mE\left[|\tilde{Y}_r|^2|\sF_t \right]\dif r\\
&&+\frac{1}{4}\mathbb{E}\left[ \int_t^T{\left\| \tilde{Z_r} \right\| ^2}\dif r\left| \mathscr{F}_t \right. \right] +\frac{1}{4}\mathbb{E}\left[ \underset{s\in \left[ t,T \right]}{\sup}\left| \tilde{Y_s} \right|^2\left| \mathscr{F}_t \right. \right] +C\mathbb{E}\left[\int_t^T{\left\| \tilde{Z_r} \right\| ^2}\dif r \left| \mathscr{F}_t \right. \right].
\de 
Combining with (\ref{iequ1}), we obtain that
\ce
&&\mathbb{E}\left[ \underset{s\in \left[ t,T \right]}{\sup}\left| \tilde{Y_s} \right|^2+\int_t^T{\left\| \tilde{Z_s} \right\|}^2\dif s\left| \mathscr{F}_t \right. \right] 
\\
&\leqslant& C\mathbb{E}\left[ \left| \tilde{\xi } \right|^2\left| \mathscr{F}_t \right. \right] +C\mathbb{E}\left[ \int_t^T|\tilde{g}\left( s,Y_{s}^{1},Z_{s}^{1} \right)|^2\dif s \left| \mathscr{F}_t \right. \right] .
\de 
The proof is complete.	    
\end{proof}

\subsection{The derivative for functions on $\cM_{2}(\mR^d)$}

In the subsection, we recall the definition of the derivative for functions on $\cM_{2}(\mR^d)$ (\cite{Lion}). 

A function $f:\cM_{2}(\mR^d)\mapsto\mR$ is differential at $\mu\in \cM_{2}(\mR^d)$, if for $\tilde{f}(\gamma):=f(\sL_\gamma),\gamma\in L^2(\Omega,\mathscr{F},\mP;\mR^d)$, there exists some $\zeta\in L^2(\Omega,\mathscr{F},\mP;\mR^d)$ with $\sL_\zeta=\mu$ such that $\tilde{f}$ is Fr\'echet differentiable at $\zeta$, that is, there exists a linear continuous mapping $D\tilde{f}(\zeta):L^2(\Omega,\mathscr{F},\mP;\mR^d)\mapsto\mR$ such that for any $\eta\in L^2(\Omega,\mathscr{F},\mP;\mR^d)$
\ce
\tilde{f}(\zeta+\eta)-\tilde{f}(\zeta)=D\tilde{f}(\zeta)(\eta)+o(|\eta|_{L^2}),  \quad|\eta|_{L^2}\rightarrow0.
\de
Since $D\tilde{f}(\zeta)\in L(L^2(\Omega,\mathscr{F},\mP;\mR^d),\mR)$, it follows from the Riesz representation theorem that there exists a $\mP$-a.s. unique variable $\vartheta\in L^2(\Omega,\mathscr{F},\mP;\mR^d)$ such that for all $\eta\in L^2(\Omega,\mathscr{F},\mP;\mR^d)$
\ce
D\tilde{f}(\zeta)(\eta)=(\vartheta,\eta)_{L^2}=\mE[\vartheta\cdot\eta].
\de

\bd\label{f1}
We say that $f\in C^1(\cM_{2}(\mR^d))$, if there exists for all $\gamma\in L^2(\Omega,\mathscr{F},\mP;\mR^d)$ a $\sL_{\gamma}$-modification of $\partial_\mu f(\sL_\gamma)(\cdot)$, again denoted by $\partial_\mu f(\sL_\gamma)(\cdot)$, such that $\partial_\mu f:\cM_{2}(\mR^d)\times\mR^d\mapsto\mR^d$ is continuous, and we identify this continuous function $\partial_\mu f$ as the derivative of $f$.
\ed

\bd\label{f2}
We say that $f\in C^2(\cM_{2}(\mR^d))$, if for any $\mu\in \cM_{2}(\mR^d)$, $f\in C^1(\cM_{2}(\mR^d))$ and $\partial_\mu f(\sL_\gamma)(\cdot)$ is differentiable, and its derivative $\partial_y\partial_\mu f:\cM_{2}(\mR^d)\times\mR^d\mapsto\mR^d\otimes\mR^d$ is continuous, and for any $y\in \mR^d$, $\partial_{\mu}f(\cdot)(y)$ is differentiable, and its derivative $\partial_{\mu}^{2}f:\cM_{2}(\mR^d)\times\mR^d\times\mR^d\mapsto\mR^d\otimes\mR^d$ is continuous.
\ed

\bd
A function $F: \mR^d\times\cM_{2}(\mR^d)\mapsto\mR$  is said to be in $C^{2,2}(\mR^d\times\cM_{2}(\mR^d))$, if

(i) $F$ is $C^2$ in $x\in\mR^d$ and $\mu\in\cM_2(\mR^d)$ respectively;

(ii) for any $\mu\in\cM_2(\mR^d)$, its derivatives 
\ce
\partial_{x}F(x,\mu), \partial_{x}^2F(x,\mu), \partial_{\mu}F(x,\mu)(y), \partial_{y}\partial_{\mu}F(x,\mu)(y), \partial_{\mu}^2F(x,\mu)(y,y^{'})
\de
are jointly continuous in the variable family $(x,\mu), (x,\mu,y)$ and $(x,\mu,y,y^{'})$ respectively.
\ed

\bd
A function $F: \mR^d\times\cM_{2}(\mR^d)\mapsto\mR$ is said to be in $C_{b}^{2,2}(\mR^d\times\cM_{2}(\mR^d))$, if
$F\in C^{2,2}(\mR^d\times\cM_{2}(\mR^d))$, and itself and all its derivatives are uniformly bounded on $\mR^d\times\cM_{2}(\mR^d)$.
\ed

\bd
A function $\Psi: [0,T]\times\mR^d\times\cM_{2}(\mR^d)\mapsto\mR$ is said to be in $C_{b}^{1,2,2}([0,T]\times\mR^d\times\cM_{2}(\mR^d))$, if

(i) $\Psi(\cdot,\cdot, \mu)\in C^{1,2}([0, T]\times\mR^d)$, for all $\mu\in\cM_{2}(\mR^d)$;

(ii) $\Psi(t,x,\cdot)\in C^2(\cM_{2}(\mR^d))$, for all $(t, x)\in[0,T]\times\mR^d$;

(iii) All derivatives of order $1$ and $2$ are continuous on $[0,T]\times\mR^d\times\cM_{2}(\mR^d)\times\mR^d$, $\partial_{\mu}\Psi$ and $\partial_y\partial_{\mu}\Psi$ are bounded over $[0,T]\times\mR^d\times\cM_{2}(\mR^d)\times\mR^d$.
\ed

\section{The existence and uniqueness of solutions to backward multivalued McKean-Vlasov SDEs}\label{teus}
In this section, we study the existence and uniqueness of  solutions for Eq.(\ref{eq1}). 

\bd\label{solu}
We say that Eq.$(\ref{eq1})$ admits a solution with the terminal value $\xi$ if there exists a triple $\left\{ \left( Y_t, K_t, Z_t\right) :t\in \left[ 0,T \right] \right\} $ which satisfies 
\begin{enumerate}[(i)]
\item 
$(Y_{\cdot},K_{\cdot})\in \sA$, $\dif \mP\times\dif t$-a.e. on $\Omega\times\left[ 0, T \right]$,
\item 
$(Y_{\cdot}, Z_{\cdot})\in \mS_{\mF}^2([0,T],\mR^d)\times\mH_{\mF}^2([0,T],\mR^{d\times l})$,
\item 
\ce
Y_t=\xi -\left( K_T-K_t \right)+\int_t^T{G\left( s, Y_s, Z_s, \sL_{(Y_s, Z_s)} \right)}\dif s-\int_t^T{Z_s}\dif W_s.
 \de
\end{enumerate}
\ed

The definition of the uniqueness for solutions to Eq.(\ref{eq1}) is the same to that for solutions to Eq.(\ref{bmsdeeq}).

In the following, we give some assumptions to assure the existence and uniqueness of  solutions for Eq.(\ref{eq1}).

\begin{enumerate}[(${\bf H}^1_G$)]
\item There exsits a non-random constant $L_G>0$ such that for a.s. $\omega\in\Omega$, it holds that
\ce
&&\left| G\left( \omega ,t,y,z,\vartheta  \right) -G\left( \omega ,t,y^{\prime},z^{\prime},\vartheta ^{\prime} \right) \right|\leq L_G\left( \left| y-y^{\prime} \right|+\| z-z^{\prime} \|+\rho \left( \vartheta ,\vartheta ^{\prime} \right) \right),\\
 &&t\in\left[ 0, T \right], y, y^{'}\in \mR^d, z, z^{'}\in \mR^{d\times l}, \vartheta, \vartheta^{'}, \in \cM_{2}(\mR^d\times\mR^{d\times l}).
\de
\end{enumerate}
\begin{enumerate}[(${\bf H}^2_G$)]
\item $\mathbb{E}\int_0^T{\left| G\left( t,0,0,\delta _{({\bf 0}, {\bf 0})} \right) \right|}^2\dif t<\infty $, where $\delta_{({\bf 0}, {\bf 0})} $ is the Dirac measure at $({\bf 0}, {\bf 0})$.
\end{enumerate}
\begin{enumerate}[(${\bf H}_{A}$)]
\item $\cD(A)=\mR^d$ and $|A^\circ(x)|\leq C(1+|x|), x\in\mR^d$, where $C>0$ is a constant.
\end{enumerate}

Now, it is the position to state the main result in the section.

\bt\label{EUni}
Assume that $({\bf H}_{A})$ holds and the coefficient $G$ satisfies $({\bf H}^1_G)$-$({\bf H}^2_G)$. Then there exists a unique solution for Eq.(\ref{eq1}).
\et

To prove the above theorem, we prepare some key lemmas. For any $\varepsilon>0$, consider the penalized backward multivalued McKean-Vlasov SDE on $\mR^d$:
\be
\left\{ \begin{array}{l}
	\mathrm{d}Y_{t}^{\varepsilon}=\[A_{\varepsilon}\left( Y_{t}^{\varepsilon}\right)-G\(t, Y_{t}^{\varepsilon}, Z_{t}^{\varepsilon}, \sL_{(Y_{t}^{\varepsilon}, Z_{t}^{\varepsilon})}\)\]\dif t+Z_{t}^{\varepsilon}\dif W_t, \quad 0\leq t<T,\\
	Y_{T}^{\varepsilon}=\xi,
\end{array}
\label{pena} 
\right. 
\ee 
where $A_{\varepsilon}$ is the Yosida approximation of  $A$. Note that $A_{\varepsilon}$ is a single-valued, maximal monotone and Lipschitz continuous function  (cf. Subsection \ref{mmo}). Thus, by Theorem \ref{bmvsdeeu}, we know that under $({\bf H}^1_G)$-$({\bf H}^2_G)$ Eq.(\ref{pena}) has a unique solution denoted as $(Y_{\cdot}^{\varepsilon}, Z_{\cdot}^{\varepsilon})$ with 
\be
\mathbb{E}\left( \underset{0\leqslant t\leqslant T}{\sup}\left| Y^\e_t \right|^2+\int_0^T{\| Z^\e_s \|^2}\dif s \right) <\infty.
\label{penaes}
\ee
Next, we give some uniform moment estimates about $Y_{\cdot}^{\varepsilon}, Z_{\cdot}^{\varepsilon}$. 

\bl\label{esti1}
Assume that ${\rm Int}(\cD(A))\neq\emptyset$ and the coefficient $G$ satisfies $({\bf H}^1_G)$-$({\bf H}^2_G)$. Then there exists a constant $C>0$ independent of $\e$ such that
\ce
\mE\left(\sup\limits_{0\leq t\leq T}|Y^\e_t|^2\right)+\mE\int_0^T\| Z^\e_s \|^2\dif s\leq C.
\de
\el
\begin{proof}
Applying the It\^o formula to $|Y_t^\e-a|^2$,  where $a$ is the same to that in Lemma \ref{yosi}, we obtain that
\be
|Y^\e_t-a|^2+\int_t^T\| Z^\e_s \|^2\dif s&=&|\xi-a|^2+2\int_t^T\left<Y^\e_s-a,G\( s, Y_{s}^{\varepsilon}, Z_{s}^{\varepsilon}, \sL_{(Y_{s}^{\varepsilon}, Z_{s}^{\varepsilon})}\)\right>\dif s\no\\
&&-2\int_t^T\<Y^\e_s-a, Z^\e_s\dif W_s\> -2\int_t^T\<Y^\e_s-a, A_{\varepsilon}\left( Y_{s}^{\varepsilon}\right)\>\dif s\no\\
&\leq&|\xi-a|^2+2\int_t^T\left<Y^\e_s-a,G\( s, Y_{s}^{\varepsilon}, Z_{s}^{\varepsilon}, \sL_{(Y_{s}^{\varepsilon}, Z_{s}^{\varepsilon})}\)\right>\dif s\no\\
&&-2\int_t^T\<Y^\e_s-a, Z^\e_s\dif W_s\>-2M_1\int_t^T|A_{\varepsilon}\left( Y_{s}^{\varepsilon}\right)|\dif s\no\\
&&+2M_2\int_t^T|Y^\e_s-a|\dif s+2M_1M_2T,
\label{itof}
\ee
where the last inequality is based on Lemma \ref{yosi}. Taking the expectation on two sides, by (\ref{penaes}), one can have that
\be
\mE|Y^\e_t-a|^2+\mE\int_t^T\| Z^\e_s \|^2\dif s&\leq&\mE|\xi-a|^2+2\mE\int_t^T|Y^\e_s-a|\left|G\( s, Y_{s}^{\varepsilon}, Z_{s}^{\varepsilon}, \sL_{(Y_{s}^{\varepsilon}, Z_{s}^{\varepsilon})}\)\right|\dif s\no\\
&&-2M_1\mE\int_t^T|A_{\varepsilon}\left( Y_{s}^{\varepsilon}\right)|\dif s+2M_2\mE\int_t^T|Y^\e_s-a|\dif s\no\\
&&+2M_1M_2T.\label{dai}
\ee
Let us estimate the second term in the right side of the above inequality. By $({\bf H}^1_G)$ it holds that
\ce
\left|G\( s, Y_{s}^{\varepsilon}, Z_{s}^{\varepsilon}, \sL_{(Y_{s}^{\varepsilon}, Z_{s}^{\varepsilon})}\)\right|&\leq&L_G\left( \left| Y_{s}^{\varepsilon}\right|+\|Z_{s}^{\varepsilon}\|+\rho \left( \sL_{(Y_{s}^{\varepsilon}, Z_{s}^{\varepsilon})},\delta _{({\bf 0}, {\bf 0})}\right) \right)+\left| G\left(s,0,0,\delta _{({\bf 0}, {\bf 0})} \right) \right|\\
&\leq&L_G\left( \left| Y_{s}^{\varepsilon}\right|+\|Z_{s}^{\varepsilon}\|+(\mE|Y_{s}^{\varepsilon}|^2)^{1/2}+(\mE\|Z_{s}^{\varepsilon}\|^2)^{1/2}\right)+\left| G\left(s,0,0,\delta _{({\bf 0}, {\bf 0})} \right) \right|\\
&\leq&L_G\left(3|a|+ \left| Y_{s}^{\varepsilon}-a\right|+\|Z_{s}^{\varepsilon}\|+(\mE|Y_{s}^{\varepsilon}-a|^2)^{1/2}+(\mE\|Z_{s}^{\varepsilon}\|^2)^{1/2}\right)\\
&&+\left| G\left(s,0,0,\delta _{({\bf 0}, {\bf 0})} \right) \right|,
\de
where we use the following fact that
\be
&&\rho^2 \left(\sL_{(Y_{s}^{\varepsilon}, Z_{s}^{\varepsilon})},\delta _{({\bf 0}, {\bf 0})} \right)\no\\
&=&\inf_{\pi\in\sC(\sL_{(Y_{s}^{\varepsilon}, Z_{s}^{\varepsilon})}, \delta _{({\bf 0}, {\bf 0})} )}\int_{(\mR^d\times\mR^{d\times l})\times(\mR^d\times\mR^{d\times l})}|(y_1,z_1)-(y_2,z_2)|^2\pi(\dif (y_1,z_1), \dif (y_2,z_2))\no\\
&\leq&\mE|(Y_{s}^{\varepsilon}, Z_{s}^{\varepsilon})-(0,0)|^2= \mE|Y_{s}^{\varepsilon}|^2+\mE\|Z_{s}^{\varepsilon}\|^2.
\label{measesti}
\ee
Thus, by  the Young inequality and the inequality $|x|\leq 1+|x|^2, x\in\mR^d$, one can obtain that
\be
&&2|Y^\e_s-a|\left|G\( s, Y_{s}^{\varepsilon}, Z_{s}^{\varepsilon}, \sL_{(Y_{s}^{\varepsilon}, Z_{s}^{\varepsilon})}\)\right|\no\\
&\leq& 6|a|L_G|Y^\e_s-a|+2L_G|Y^\e_s-a|^2+2L_G|Y^\e_s-a|\|Z_{s}^{\varepsilon}\|+2L_G|Y^\e_s-a|(\mE|Y_{s}^{\varepsilon}-a|^2)^{1/2}\no\\
&&+2L_G|Y^\e_s-a|(\mE\|Z_{s}^{\varepsilon}\|^2)^{1/2}+2|Y^\e_s-a|\left| G\left(s,0,0,\delta _{({\bf 0}, {\bf 0})} \right) \right|\no\\
&\leq& 6|a|L_G+(6|a|L_G+2L_G+9L^2_G+1)|Y^\e_s-a|^2+\mE|Y_{s}^{\varepsilon}-a|^2\no\\
&&+\frac{1}{4}\|Z_{s}^{\varepsilon}\|^2+\frac{1}{4}\mE\|Z_{s}^{\varepsilon}\|^2+\left| G\left(s,0,0,\delta _{({\bf 0}, {\bf 0})} \right) \right|^2.\label{tiao}
\ee

Combining (\ref{dai}) with (\ref{tiao}), by $({\bf H}^2_G)$ we get that
\be
\mE|Y^\e_t-a|^2+\frac{1}{2}\mE\int_t^T\| Z^\e_s \|^2\dif s\leq C+C\int_t^T\mE|Y^\e_s-a|^2\dif s,
\label{zees}
\ee
where $C$ is independent of $\e$. By the Gronwall inequality, it holds that 
\ce
\sup\limits_{0\leq t\leq T}\mE|Y^\e_t-a|^2\leq Ce^{CT}.
\de
From this and (\ref{zees}), it follows that 
\be
\mE\int_0^T\| Z^\e_s \|^2\dif s\leq C.
\label{zest}
\ee

Next, for (\ref{itof}), by the BDG inequality and the Young inequality we have that
\ce
\mE\left(\sup\limits_{0\leq t\leq T}|Y^\e_t-a|^2\right)&\leq&\mE|\xi-a|^2+2\mE\int_0^T\left|\<Y^\e_s-a,G\( s, Y_{s}^{\varepsilon}, Z_{s}^{\varepsilon}, \sL_{(Y_{s}^{\varepsilon}, Z_{s}^{\varepsilon})}\)\>\right|\dif s\\
&&+2\mE\left(\sup\limits_{0\leq t\leq T}\left|\int_t^T\<Y^\e_s-a, Z^\e_s\dif W_s\>\right|\right)+2M_2\mE\int_0^T|Y^\e_s-a|\dif s\\
&&+2M_1M_2T\\
&\leq&\mE|\xi-a|^2+2\mE\int_0^T\left|\<Y^\e_s-a,G\( s, Y_{s}^{\varepsilon}, Z_{s}^{\varepsilon}, \sL_{(Y_{s}^{\varepsilon}, Z_{s}^{\varepsilon})}\)\>\right|\dif s\\
&&+2C\mE\left(\int_0^T|Y^\e_s-a|^2\|Z^\e_s\|^2\dif s \right)^{1/2}+2M_2\mE\int_0^T|Y^\e_s-a|^2\dif s\\
&&+2(M_1+1)M_2T\\
&\leq&\mE|\xi-a|^2+2\mE\int_0^T\left|\<Y^\e_s-a,G\( s, Y_{s}^{\varepsilon}, Z_{s}^{\varepsilon}, \sL_{(Y_{s}^{\varepsilon}, Z_{s}^{\varepsilon})}\)\>\right|\dif s\\
&&+\frac{1}{2}\mE\left(\sup\limits_{0\leq t\leq T}|Y^\e_t-a|^2\right)+C\mE\int_0^T\| Z^\e_s \|^2\dif s+2M_2\mE\int_0^T|Y^\e_s-a|^2\dif s\\
&&+2(M_1+1)M_2T.
\de
From (\ref{tiao}) and (\ref{zest}), it follows that
\be
\mE\left(\sup\limits_{0\leq t\leq T}|Y^\e_t|^2\right)\leq 2\mE\left(\sup\limits_{0\leq t\leq T}|Y^\e_t-a|^2\right)+2|a|^2\leq C.
\label{yest}
\ee
The proof is complete.
\end{proof}

\bl\label{esti2}
Assume that $({\bf H}_{A})$ holds and the coefficient $G$ satisfies $({\bf H}^1_G)$-$({\bf H}^2_G)$. Then there exists a constant $C>0$ independent of $\e$ such that
\ce
\mE\int_0^T|A_{\varepsilon}\left( Y_{s}^{\varepsilon}\right)|\dif s+\mE\int_0^T|A_{\varepsilon}\left( Y_{s}^{\varepsilon}\right)|^2\dif s\leq C.
\de
\el
\begin{proof}
By (\ref{itof}), it holds that
\ce
2M_1\int_t^T|A_{\varepsilon}\left( Y_{s}^{\varepsilon}\right)|\dif s&\leq&|\xi-a|^2+2\int_t^T\left<Y^\e_s-a,G\( s, Y_{s}^{\varepsilon}, Z_{s}^{\varepsilon}, \sL_{(Y_{s}^{\varepsilon}, Z_{s}^{\varepsilon})}\)\right>\dif s\no\\
&&-2\int_t^T\<Y^\e_s-a, Z^\e_s\dif W_s\>+2M_2\int_t^T|Y^\e_s-a|\dif s+2M_1M_2T,
\de
and furthermore
\ce
2M_1\mE\int_t^T|A_{\varepsilon}\left( Y_{s}^{\varepsilon}\right)|\dif s&\leq&\mE|\xi-a|^2+2\mE\int_t^T\left<Y^\e_s-a,G\( s, Y_{s}^{\varepsilon}, Z_{s}^{\varepsilon}, \sL_{(Y_{s}^{\varepsilon}, Z_{s}^{\varepsilon})}\)\right>\dif s\no\\
&&+2M_2\mE\int_t^T|Y^\e_s-a|\dif s+2M_1M_2T,
\de
where (\ref{penaes}) is used. Thus, by (\ref{tiao}), (\ref{zest}) and (\ref{yest}), we know that
\ce
\mE\int_t^T|A_{\varepsilon}\left( Y_{s}^{\varepsilon}\right)|\dif s\leq C,
\de
where the constant $C$ is independent of $\e$. 

Next, note that $|A_{\varepsilon}(x)|\leq |A^\circ(x)|$ for $x\in \cD(A)$. Thus, by $({\bf H}_A)$ and (\ref{yest}), it holds that
\ce
\mE\int_0^T|A_{\varepsilon}\left(Y_{s}^{\varepsilon}\right)|^2\dif s\leq \mE\int_0^T|A^\circ\left(Y_{s}^{\varepsilon}\right)|^2\dif s\leq C\mE\int_0^T(1+|Y_{s}^{\varepsilon}|)^2\dif s\leq C.
\de
The proof is complete.
\end{proof}

\bl\label{esti3}
Assume that $({\bf H}_{A})$ holds and the coefficient $G$ satisfies $({\bf H}^1_G)$-$({\bf H}^2_G)$. Then for any $\e, \d>0$, there exists a constant $C>0$ independent of $\e, \d$ such that
\ce
\mE\left(\sup\limits_{0\leq t\leq T}|Y^\e_t-Y^\d_t|^2\right)+\mE\int_0^T\| Z^\e_s-Z^\d_s \|^2\dif s\leq C(\e+\d).
\de
\el
\begin{proof}
For any $\e, \d>0$, it holds that
\ce
Y^\e_t-Y^\d_t&=&\int_t^T\[G(s, Y_s^\e, Z_s^\e, \sL_{(Y_s^\e, Z_s^\e)})-G(s, Y_s^\d, Z_s^\d, \sL_{(Y_s^\d, Z_s^\d)})\]\dif s\\
&&-\int_t^T\(A_\e(Y_s^\e)-A_\d(Y_s^\d)\)\dif s-\int_t^T(Z_s^\e-Z_s^\d)\dif W_s.
\de
Applying the It\^o formula to $|Y^\e_t-Y^\d_t|^2$, we know that
\be
&&|Y^\e_t-Y^\d_t|^2+\int_t^T\|Z^\e_s-Z^\d_s\|^2\dif s\no\\
&=&2\int_t^T\<Y^\e_s-Y^\d_s,G(s, Y_s^\e, Z_s^\e, \sL_{(Y_s^\e, Z_s^\e)})-G(s, Y_s^\d, Z_s^\d, \sL_{(Y_s^\d, Z_s^\d)})\>\dif s\no\\
&&-2\int_t^T\<Y^\e_s-Y^\d_s,A_\e(Y_s^\e)-A_\d(Y_s^\d)\>\dif s-2\int_t^T\<Y^\e_s-Y^\d_s,(Z_s^\e-Z_s^\d)\dif W_s\>,
\label{itofed}
\ee
and
\be
&&\mE|Y^\e_t-Y^\d_t|^2+\mE\int_t^T\|Z^\e_s-Z^\d_s\|^2\dif s\no\\
&=&2\mE\int_t^T\<Y^\e_s-Y^\d_s,G(s, Y_s^\e, Z_s^\e, \sL_{(Y_s^\e, Z_s^\e)})-G(s, Y_s^\d, Z_s^\d, \sL_{(Y_s^\d, Z_s^\d)})\>\dif s\no\\
&&-2\mE\int_t^T\<Y^\e_s-Y^\d_s,A_\e(Y_s^\e)-A_\d(Y_s^\d)\>\dif s.
\label{itofedex}
\ee

From $({\bf H}^1_G)$ and (\ref{measesti}), it follows that
\be
&&2|\<Y^\e_s-Y^\d_s,G(s, Y_s^\e, Z_s^\e, \sL_{(Y_s^\e, Z_s^\e)})-G(s, Y_s^\d, Z_s^\d, \sL_{(Y_s^\d, Z_s^\d)})\>|\no\\
&\leq& 2|Y^\e_s-Y^\d_s||G(s, Y_s^\e, Z_s^\e, \sL_{(Y_s^\e, Z_s^\e)})-G(s, Y_s^\d, Z_s^\d, \sL_{(Y_s^\d, Z_s^\d)})|\no\\
&\leq&2 |Y^\e_s-Y^\d_s|L_G(|Y_s^\e-Y_s^\d|+\|Z_s^\e-Z_s^\d\|+\rho(\sL_{(Y_s^\e, Z_s^\e)},\sL_{(Y_s^\d, Z_s^\d)}))\no\\
&\leq&2 |Y^\e_s-Y^\d_s|L_G(|Y_s^\e-Y_s^\d|+\|Z_s^\e-Z_s^\d\|+(\mE|Y_s^\e-Y_s^\d|^2)^{1/2}+(\mE\|Z_s^\e-Z_s^\d\|^2)^{1/2})\no\\
&\leq&(2L_G+9L_G^2) |Y^\e_s-Y^\d_s|^2+\frac{1}{4}\|Z_s^\e-Z_s^\d\|^2+\mE|Y_s^\e-Y_s^\d|^2\no\\
&&+\frac{1}{4}\mE\|Z_s^\e-Z_s^\d\|^2.
\label{gtiao}
\ee
Moreover, by the deduction in \cite[Proposition 6]{MODE}, we know that 
\be
-2\<Y^\e_s-Y^\d_s,A_\e(Y_s^\e)-A_\d(Y_s^\d)\>\leq 3\e|A_\e(Y_s^\e)|^2+3\d|A_\d(Y_s^\d)|^2.
\label{yeydyosi}
\ee
Inserting (\ref{gtiao}) and (\ref{yeydyosi}) in (\ref{itofedex}), by Lemma \ref{esti2} one can obtain that
\ce
\mE|Y^\e_t-Y^\d_t|^2+\frac{1}{2}\mE\int_t^T\|Z^\e_s-Z^\d_s\|^2\dif s&\leq& (2L_G+9L_G^2+1)\int_t^T\mE |Y^\e_s-Y^\d_s|^2\dif s\\
&&+\mE\int_t^T(3\e|A_\e(Y_s^\e)|^2+3\d|A_\d(Y_s^\d)|^2)\dif s\\
&\leq& (2L_G+9L_G^2+1)\int_t^T\mE |Y^\e_s-Y^\d_s|^2\dif s+C(\e+\d).
\de
Thus, by the Gronwall inequality, it holds that
\be
\sup\limits_{0\leq t\leq T}\mE|Y^\e_t-Y^\d_t|^2\leq C(\e+\d),
\label{supyeyd}
\ee
and furthermore
\be
\mE\int_0^T\|Z^\e_s-Z^\d_s\|^2\dif s\leq C(\e+\d).
\label{zezd}
\ee

Next, from (\ref{itofed}) and (\ref{gtiao})-(\ref{zezd}), it follows that
\ce
&&\mE\left(\sup\limits_{0\leq t\leq T}|Y^\e_t-Y^\d_t|^2\right)\\
&\leq& 2\mE\int_0^T|\<Y^\e_s-Y^\d_s,G(s, Y_s^\e, Z_s^\e, \sL_{(Y_s^\e, Z_s^\e)})-G(s, Y_s^\d, Z_s^\d, \sL_{(Y_s^\d, Z_s^\d)})\>|\dif s\\
&&+\mE\int_0^T(3\e|A_\e(Y_s^\e)|^2+3\d|A_\d(Y_s^\d)|^2)\dif s+2\mE\left(\sup\limits_{0\leq t\leq T}\left|\int_t^T\<Y^\e_s-Y^\d_s,(Z_s^\e-Z_s^\d)\dif W_s\>\right|\right)\\
&\leq&(2L_G+9L_G^2+1)\int_0^T\mE |Y^\e_s-Y^\d_s|^2\dif s+\frac{1}{2}\mE\int_t^T\|Z^\e_s-Z^\d_s\|^2\dif s+C(\e+\d)\\
&&+C\mE\left(\int_0^T |Y^\e_s-Y^\d_s|^2\|Z^\e_s-Z^\d_s\|^2\dif s\right)^{1/2}\\
&\leq&C(\e+\d)+\frac{1}{2}\mE\left(\sup\limits_{0\leq t\leq T}|Y^\e_t-Y^\d_t|^2\right)+C\mE\int_0^T\|Z^\e_s-Z^\d_s\|^2\dif s,
\de
which yields that
\ce
\mE\left(\sup\limits_{0\leq t\leq T}|Y^\e_t-Y^\d_t|^2\right)\leq C(\e+\d).
\de
The proof is complete.
\end{proof}

\medspace

Now, we apply the above lemmas to prove Theorem \ref{EUni}.

{\bf Proof of the existence for Theorem \ref{EUni}.}

We prove that the limit $(Y_{\cdot},K_{\cdot},Z_{\cdot})$ of the sequence $\{(Y_{\cdot}^{\e},K_{\cdot}^{\e},Z_{\cdot}^{\e})\}$ is a solution of Eq.(\ref{eq1}), where $K_t^\e:=\int_0^t A_\e(Y_s^\e)\dif s$.

From Lemma \ref{esti3}, it follows that $\left\lbrace Y_{\cdot}^{\e} \right\rbrace $ is a Cauchy sequence in $L^2\left(  \Omega, \sF, \mP; C([0,T], \mR^d) \right)$ and $\left\lbrace Z_{\cdot}^{\e} \right\rbrace $ is a Cauchy sequence in $L^2\left( \Omega\times\left[ 0,T \right], \dif \mathbb{P}\times \dif t \right)$. So, there exist two processes $Y, Z$ such that  
\be
&&\underset{\e\rightarrow 0}{\lim}\mathbb{E}\left(\sup\limits_{0\leq t\leq T}|Y^\e_t-Y_t|^2\right) =0, \label{equ6}\\
&&\underset{\e\rightarrow 0}{\lim}\mathbb{E}\int_0^T{\left| Z_{s}^{\e}-Z_{s} \right|}^2\dif s=0.\label{equ7}
\ee
Moreover, by (\ref{measesti}) it holds that
\be
\underset{\e\rightarrow 0}{\lim}\rho\left( \sL_{(Y_{t}^{\e},Z_{t}^{\e})},\sL_{(Y_{t}, Z_t)} \right) \leqslant \underset{\e\rightarrow 0}{\lim}\((\mathbb{E}|Y^\e_t-Y_t|^2)^{1/2}+(\mathbb{E}\|Z^\e_t-Z_t\|^2)^{1/2}\)=0. \label{equ8}
\ee

Now, let us put 
\ce
-K_t:=Y_0-Y_t-\int_0^t{G\left( s, Y_s, Z_s, \sL_{(Y_{s}, Z_s)}\right)}\dif s+\int_0^t{Z_s}\dif W_s.
\de
(\ref{equ6})-(\ref{equ8}) imply that 
\ce
\underset{\e\rightarrow 0}{\lim}\mathbb{E}\left( \underset{t\in \left[ 0,T \right]}{\mathrm{sup}}\left| K_{t}^{\e}-K_t \right|^2 \right) =0,
\de
which together with Lemma \ref{limiconv}, gives that
\be
\underset{\e\rightarrow 0}{\lim}\int_0^t{\left< Y_{s}^{\e},\dif K_{s}^{\e} \right>}=\int_0^t{\left< Y_s,\dif K_s \right>}.\label{equ10}
\ee
Note that $(Y_{\cdot}^{\e},K_{\cdot}^{\e})\in\sA$. Thus, by Lemma \ref{equi}, it holds that for any $(Y^{'},K^{'})\in\sA$,
$$
\left\langle Y_{t}^{\e}-Y_{t}^{'},\dif K_{t}^{\e}-\dif K_{t}^{'}\right\rangle \geq0,
$$
which together with (\ref{equ10}), yields $(Y_{\cdot},K_{\cdot})\in\sA, a.s. $. That is, $\left( Y_t,K_t,Z_t \right) _{t\in \left[ 0,T \right]}$ is a solution of Eq.(\ref{eq1}).

\medspace

To show the uniqueness part in Theorem \ref{EUni}, we prepare the following proposition.

\bp\label{diffesti2}
Assume that $({\bf H}_{A})$ holds. Let $G^1, G^2$ be two functions satisfying the standard assumptions $({\bf H}^1_G)$-$({\bf H}^2_G)$. Let $\xi^1, \xi^2 \in L^{2}(\Omega, \sF_T, \mP)$ and $(Y^i, K^i, Z^i), i=1,2,$ be the solutions of Eq.(\ref{eq1}) with the terminal values $\xi^1, \xi^2$, respectively. Set $\overline{Y}_t:=Y^1_t-Y^2_t, \overline{Z}_t:=Z^1_t-Z^2_t, \overline{K}_t:=K^1_t-K^2_t, \overline{\xi}:=\xi ^1-\xi ^2, \bar{G}:=G^1-G^2$. Then it holds that
\be
&&\mE\int_0^T(|\overline{Y}_t |^2+\|\bar{Z}_t\|^2)\dif t\leq C\Gamma,\label{barest1}\\
&&\mE\(\sup\limits_{t\in[0,T]}|\overline{Y}_t |^2\)\leq C\Gamma,\label{barest2}
\ee
where 
$$
\Gamma:=\mE|\overline{\xi }|^2+\mathbb{E}\int_0^T|\bar{G}( s,Y_{s}^{1},Z_{s}^{1},\sL_{(Y^1_s,Z^1_s)})|^2\dif s.
$$
\ep
\begin{proof}
Since $(Y^i, K^i, Z^i), i=1,2,$ are the solutions of Eq.(\ref{eq1}) with the terminal values $\xi^1, \xi^2$, respectively, it holds that
\ce
&&Y^1_t=\xi^1 -\left( K^1_T-K^1_t \right)+\int_t^T{G^1\left( s,Y^1_s,Z^1_s,\sL_{(Y^1_s,Z^1_s)}\right)}\dif s-\int_t^T{Z^1_s}\dif W_s,\\
&&Y^2_t=\xi^2 -\left( K^2_T-K^2_t \right)+\int_t^T{G^2\left( s,Y^2_s,Z^2_s,\sL_{(Y^2_s,Z^2_s)}\right)}\dif s-\int_t^T{Z^2_s}\dif W_s,
\de
and 
\ce
\overline{Y}_t&=&\overline{\xi }-\int_t^T{\dif \overline{K}_s}+\int_t^T\(G^1\left( s,Y_{s}^{1},Z_{s}^{1},\sL_{(Y^1_s,Z^1_s)}\right)-G^2\left( s,Y_{s}^{2},Z_{s}^{2},\sL_{(Y^2_s,Z^2_s)}\right)\) \dif s
\\
&&-\int_t^T{\overline{Z}_s\dif W_s}.
\de 	
Applying the It\^{o} formula to $|\overline{Y}_t|^2$, by Lemma \ref{equi} one can obtain that
\be
&&|\overline{Y}_t |^2+\int_t^T\|\bar{Z}_s\|^2\dif s\no\\
&=&|\overline{\xi }|^2+2\int_t^T\<\overline{Y}_s, G^1( s,Y_{s}^{1},Z_{s}^{1},\sL_{(Y^1_s,Z^1_s)})-G^2( s,Y_{s}^{2},Z_{s}^{2},\sL_{(Y^2_s,Z^2_s)})\>\dif s\no\\
&&-2\int_t^T\<\overline{Y}_s, \overline{Z}_s\dif W_s\>-2\int_t^T\<\overline{Y}_s, \dif \overline{K}_s\>\no\\
&\leq&|\overline{\xi }|^2+2\int_t^T\<\overline{Y}_s, \bar{G}( s,Y_{s}^{1},Z_{s}^{1},\sL_{(Y^1_s,Z^1_s)})\>\dif s-2\int_t^T\<\overline{Y}_s, \overline{Z}_s\dif W_s\>\no\\
&&+2\int_t^T\<\overline{Y}_s, G^2( s,Y_{s}^{1},Z_{s}^{1},\sL_{(Y^1_s,Z^1_s)})-G^2( s,Y_{s}^{2},Z_{s}^{2},\sL_{(Y^2_s,Z^2_s)})\>\dif s.
\label{barest3}
\ee 	
Taking the expectation on two sides of the above inequality, we get that
\be
&&\mE|\overline{Y}_t |^2+\mE\int_t^T\| \bar{Z}_s\|^2\dif s\no\\
&\leq&\mE|\overline{\xi }|^2+2\mE\int_t^T\<\overline{Y}_s, \bar{G}( s,Y_{s}^{1},Z_{s}^{1},\sL_{(Y^1_s,Z^1_s)})\>\dif s\no\\
&&+2\mE\int_t^T\<\overline{Y}_s, G^2( s,Y_{s}^{1},Z_{s}^{1},\sL_{(Y^1_s,Z^1_s)})-G^2( s,Y_{s}^{2},Z_{s}^{2},\sL_{(Y^2_s,Z^2_s)})\>\dif s\no\\
&\leq& \mE|\overline{\xi }|^2+\int_t^T\mE|\overline{Y}_s|^2\dif s+\int_t^T\mE|\bar{G}( s,Y_{s}^{1},Z_{s}^{1},\sL_{(Y^1_s,Z^1_s)})|^2\dif s\no\\
&&+2\mathbb{E}\int_t^T{|\overline{Y}_s|\cdot L_{G^2}\( \left| \overline{Y}_{s}\right|+\| \bar{Z}_{s} \|+\rho \left( \sL_{(Y^1_s,Z^1_s)},\sL_{(Y^2_s,Z^2_s)} \right) \)}\dif s\no\\
&=& \mE|\overline{\xi }|^2+\int_t^T\mE|\overline{Y}_s|^2\dif s+\int_t^T\mE|\bar{G}( s,Y_{s}^{1},Z_{s}^{1},\sL_{(Y^1_s,Z^1_s)})|^2\dif s\no\\
&&+2L_{G^2}\mathbb{E}\int_t^T{\left| \overline{Y}_s\right|^2}\dif s+2L_{G^2}\mathbb{E}\int_t^T{\left| \overline{Y}_s\right|\| \bar{Z}_s\|}\dif s\no\\
&&+2L_{G^2}\mathbb{E}\int_t^T{ \left| \overline{Y}_s\right|\rho \left( \sL_{(Y^1_s,Z^1_s)},\sL_{(Y^2_s,Z^2_s)} \right) }\dif s\no\\
&\leqslant& \mE|\overline{\xi }|^2+\int_t^T\mE|\overline{Y}_s|^2\dif s+\int_t^T\mE|\bar{G}( s,Y_{s}^{1},Z_{s}^{1},\sL_{(Y^1_s,Z^1_s)})|^2\dif s\no\\ 
&&+C\mathbb{E}\int_t^T{\left| \overline{Y}_s \right|^2}\dif s+\frac{1}{4}\mathbb{E}\int_t^T{\| \bar{Z}_s \|^2}\dif s\no\\
&&+2L_{G^2}\mathbb{E}\int_t^T\left| \overline{Y}_s\right|\left(\left(\mE\left| \overline{Y}_s \right|^2\right)^{1/2}+\left(\mE\| \bar{Z}_s \|^2\right)^{1/2}\right)\dif s\no\\
&\leq&\mE|\overline{\xi }|^2+\int_t^T\mE|\overline{Y}_s|^2\dif s+\int_t^T\mE|\bar{G}( s,Y_{s}^{1},Z_{s}^{1},\sL_{(Y^1_s,Z^1_s)})|^2\dif s\no\\
&&+C\int_t^T\mE|\overline{Y}_s|^2\dif s+\frac{1}{2}\mathbb{E}\int_t^T\|\bar{Z}_s\|^2\dif s,
\label{barest4}
\ee
and
\ce
\mE|\overline{Y}_t |^2\leq \mE|\overline{\xi }|^2+\int_t^T\mE|\overline{Y}_s|^2\dif s+\int_t^T\mE|\bar{G}( s,Y_{s}^{1},Z_{s}^{1},\sL_{(Y^1_s,Z^1_s)})|^2\dif s,
\de
which together with the Gronwall inequality yields that
\be
\sup\limits_{t\in[0,T]}\mE|\overline{Y}_t |^2\leq C\left(\mE|\overline{\xi }|^2+\int_0^T\mE|\bar{G}( s,Y_{s}^{1},Z_{s}^{1},\sL_{(Y^1_s,Z^1_s)})|^2\dif s\right).
\label{barest5}
\ee
Inserting (\ref{barest5}) in (\ref{barest4}), we have (\ref{barest1}).

Next, from (\ref{barest3}) and the BDG inequality, it follows that
\ce
&&\mE(\sup\limits_{t\in[0,T]}|\overline{Y}_t |^2)\\
&\leq& \mE|\overline{\xi }|^2+2\int_0^T|\<\overline{Y}_s, \bar{G}( s,Y_{s}^{1},Z_{s}^{1},\sL_{(Y^1_s,Z^1_s)})\>|\dif s+2\mE\left(\sup\limits_{t\in[0,T]}|\int_t^T\<\overline{Y}_s, \overline{Z}_s\dif W_s\>|\right)\\
&&+2\mE\int_0^T|\<\overline{Y}_s, G^2( s,Y_{s}^{1},Z_{s}^{1},\sL_{(Y^1_s,Z^1_s)})-G^2( s,Y_{s}^{2},Z_{s}^{2},\sL_{(Y^2_s,Z^2_s)})\>|\dif s\\
&\leq& \mE|\overline{\xi }|^2+\mE\int_0^T|\overline{Y}_s|^2\dif s+\int_0^T\mE|\bar{G}( s,Y_{s}^{1},Z_{s}^{1},\sL_{(Y^1_s,Z^1_s)})|^2\dif s+C\mE\left(\int_0^T|\overline{Y}_t|^2\| \bar{Z}_t\|^2\dif t\right)^{1/2}\\
&&+C\int_0^T\mE|\overline{Y}_s|^2\dif s+\frac{1}{2}\mathbb{E}\int_0^T\|\bar{Z}_s\|^2\dif s\\
&\leq& \mE|\overline{\xi }|^2+\mE\int_0^T|\overline{Y}_s|^2\dif s+\int_0^T\mE|\bar{G}( s,Y_{s}^{1},Z_{s}^{1},\sL_{(Y^1_s,Z^1_s)})|^2\dif s+\frac{1}{2}\mE\left(\sup\limits_{t\in[0,T]}|\overline{Y}_t |^2\right)\\
&&+C\mE\int_0^T\| \bar{Z}_s\|^2\dif s+C\int_0^T\mE|\overline{Y}_s|^2\dif s+\frac{1}{2}\mathbb{E}\int_0^T\|\bar{Z}_s\|^2\dif s,
\de
which together with (\ref{barest1}) gives (\ref{barest2}). The proof is complete.
\end{proof}

{\bf The proof of the uniqueness for Theorem \ref{EUni}.} 

Assume that $\left(Y^1_{\cdot}, K^1_{\cdot}, Z^1_{\cdot} \right) $ and $\left(Y^2_{\cdot}, K^2_{\cdot}, Z^2_{\cdot} \right)$ are two  solutions of Eq.(\ref{eq1}) with $Y^1_T=Y^2_T=\xi, K_T=\tilde{K}_T$, i.e.
\ce
&&Y^1_t=\xi -\left( K^1_T-K^1_t \right)+\int_t^T{G\left( s,Y^1_s,Z^1_s,\sL_{(Y^1_s, Z^1_s)}\right)}\dif s-\int_t^T{Z^1_s}\dif W_s,\\
&&Y^2_t=\xi -\left( K^2_T-K^2_t \right)+\int_t^T{G\left( s,Y^2_s,Z^2_s,\sL_{(Y^2_s, Z^2_s)}\right)}\dif s-\int_t^T{Z^2_s}\dif W_s.
\de
By Proposition \ref{diffesti2} with $\xi^1=\xi^2, G^1=G^2$, it holds that
\ce
\mathbb{E}\left(\sup\limits_{t\in[0,T]}\left| \bar{Y}_t\right|^2\right)=0, \,\, \mathbb{E}\int_0^T{\| \bar{Z}_s\|}^2\dif s = 0,
\de
which gives that 
\ce
Y^1_{t}=Y^2_{t}, \quad t\in[0,T], \,\, \mP-a.s., \quad Z^1_{\cdot}=Z^2_{\cdot},\dif \mathbb{P}\times dt-a.s..
\de

Finally, from the above deduction it follows that for any $t\in[0,T]$,
\ce
-K^1_t&=&\xi- Y^1_t +\int_t^T{G\left( s,Y^1_s,Z^1_s, \sL_{(Y^1_{s}, Z^1_s)} \right)}\dif s-\int_t^T{Z^1_s}\dif W_s- K^1_T\\
&=&\xi- Y^2_t +\int_t^T{G\left( s,Y^2_s,Z^2_s, \sL_{(Y^2_{s}, Z^2_s)} \right)}\dif s-\int_t^T{Z^2_s}\dif W_s- K^2_T \\
&=&-K^2_t.
\de
So, the fact that $K_t$ is continuous in $t$ assures that $K^1_t=K^2_t, t\in[0,T] \,\, \mP$-a.s.. The proof is complete.

\medspace

By Proposition \ref{diffesti2}, we also have the following result.

\bc\label{continit}
Assume that $({\bf H}_{A})$ holds and $G$ satisfies $({\bf H}^1_G)$-$({\bf H}^2_G)$. Let $\xi^1, \xi^2 \in L^{2}(\Omega, \sF_T, \mP)$ and $(Y^i, K^i, Z^i), i=1,2,$ be the solutions of Eq.(\ref{eq1}) with the terminal values $\xi^1, \xi^2$, respectively. Then it holds that
\ce
\mE\(\sup\limits_{t\in[0,T]}|\overline{Y}_t |^2\)\leq C\mE|\overline{\xi }|^2.
\de
\ec

By the above corollary, we know that the solution of Eq.(\ref{eq1}) depends continuously on the terminal value.

\section{Connection with parabolic variational inequalities}\label{CONN}

In this section, we apply the result in the previous section to a type of parabolic variational inequalities and give a probabilistic representation of their solutions.

First of all, we construct a new filtration. Let $\{\tilde{\sF}_t\}_{t\in[0,T]}$ be the filtration generated by $(W_t)_{t\in[0,T]}$ and augmented by a $\sigma$-field $\sF^0$, i.e.,
\ce
\sF^W_t:=\sigma\{W_s: 0\leq s\leq t\}, \quad \tilde{\sF}_t:=\(\bigcap\limits_{s>t}\sF^W_s\)\vee\sF^0, \quad t\in[0,T],
\de
where $\sF^0\subset\sF$ has the following properties:

(i) $(W_t)_{t\in[0,T]}$ is independent of $\sF^0$;

(ii) $\cM_2(\mR^m)=\{\mP\circ \eta^{-1}, \eta\in L^2(\sF^0; \mR^m)\}$;

(iii) $\sF^0\supset \cN$ and $\cN$ is the collection of all $\mP$-null sets.

In addition, consider the following forward McKean-Vlasov SDEs on $\mR^m$: for any $t\in[0,T]$
\be \left \{ \begin{array}{l}
\dif X_{s}^{t,\eta}= b\left(s,X_{s}^{t,\eta},\sL_{X _{s}^{t,\eta}} \right)\dif s+ \sigma \left(s,X_{s}^{t,\eta},\sL_{X _{s}^{t,\eta}} \right)\dif W_s, t\leq s\leq T, \hfill \\
X_{t}^{t,\eta}=\eta,
\end{array}
\label{equ11}	
\right .
\ee
and
\be \left \{ \begin{array}{l}
\dif X_s^{t,x,\eta}= b(s,X^{t,x,\eta}_s,\sL_{X^{t,\eta}_s})\dif s+\sigma(s,X^{t,x,\eta}_s,\sL_{X^{t,\eta}_s})\dif W_s,t\leq s\leq T, \hfill \\
X_{t}^{t,x,\eta}=x,
\end{array}
\label{equ111}	
\right .
\ee
where $\eta$ is a $\tilde{\sF}_t$-measurable random variable with $\mE|\eta|^2<\infty$, and the coefficients $b:\left[ 0,T\right]\times \mR^m\times\cM_2(\mR^m)\mapsto{\mR^m}$, $\sigma:\left[ 0,T\right]\times\mR^m\times\cM_2(\mR^m)\mapsto{\mR^{m\times l}}$ are Borel measurable. 

Assume:
\begin{enumerate}[(${\bf H}^1_{b,\sigma}$)] 
\item The functions $b, \sigma$ are continuous in $(t,x,\mu)$ and satisfy for $(t,x,\mu)\in\left[ 0,T\right]\times\mR^{m}\times{\cM_{2}(\mR^m)}$
\ce
{|{b(t,x,\mu)}|}+\|\sigma(t,x,\mu)\|\leq{L_{b,\sigma}(t)(|{x}|+\|{\mu}\|_{2})},
\de
where $L_{b,\sigma}: [0, T]\mapsto (0, \infty)$ is an increasing function.
\end{enumerate}
\begin{enumerate}[(${\bf H}^2_{b,\sigma}$)] 
	\item The functions $b,\sigma$ satisfy for $t\in[0,T], (x_1,\mu_1), (x_2,\mu_2)\in\mR^{m}\times{\cM_{2}(\mR^m)}$
	\ce
	\left| b(t,x_1,\mu_1)-b(t,x_2,\mu_2)\right| +\parallel{\sigma(t,x_1,\mu_1)-\sigma(t,x_2,\mu_2)}\parallel
	\leq{L_{b,\sigma}(t)\(|x_1-x_2|+\rho(\mu_1,\mu_2)\)}.
	\de
\end{enumerate}

By Theorem 3.1 in \cite{dq1}, it holds that under  $({\bf H}^1_{b,\sigma})$-$({\bf H}^2_{b,\sigma})$, Eq.(\ref{equ11}) has a unique solution denoted as $X_{\cdot}^{t,\eta}$. Thus, by inserting $\sL_{X_{\cdot}^{t,\eta}}$ in Eq.(\ref{equ111}), it becomes a classical SDE. Based on \cite[Theorem 19.3]{hzy}, we know that under $({\bf H}^1_{b,\sigma})$-$({\bf H}^2_{b,\sigma})$, Eq.(\ref{equ111}) has a unique solution denoted as $X_{\cdot}^{t,x,\eta}$. Moreover, about $X_{\cdot}^{t,\eta}, X_{\cdot}^{t,x,\eta}$ we have the following result (\cite[Lemma 3.1]{L}). 

\bl\label{xest}
Suppose that $({\bf H}^1_{b,\sigma})$-$({\bf H}^2_{b,\sigma})$ hold. Then there exists a constant $C>0$ such that for any $t\in[0,T]$, $x,\bar{x}\in\mR^m$, $\eta, \bar{\eta}\in L^2(\tilde{\sF}_t,\mR^m)$,

(i) $\mE\left[\sup\limits_{s\in[t,T]}|X_s^{t,\eta}-X_s^{t,\bar{\eta}}|^2|\tilde{\sF}_t\right]\leq C\left(|\eta-\bar{\eta}|^2+\rho^2(\sL_\eta,\sL_{\bar{\eta}})\right)$,

(ii) $\mE\left[\sup\limits_{s\in[t,T]}|X_s^{t,x,\eta}-X_s^{t,\bar{x},\bar{\eta}}|^2|\tilde{\sF}_t\right]\leq C\left(|x-\bar{x}|^2+\rho^2(\sL_\eta,\sL_{\bar{\eta}})\right)$,

(iii) $\sup\limits_{s\in[t,T]}\rho(\sL_{X_s^{t,\eta}},\sL_{X_s^{t,\bar{\eta}}})\leq C\rho(\sL_\eta,\sL_{\bar{\eta}})$.
\el

Next, consider the following backward multivalued McKean-Vlasov SDEs on $\mR$:
\be \left \{ \begin{array}{l}
\dif Y_{s}^{t,\eta} \in A'\left( Y_{s}^{t,\eta} \right) \dif s-H\left( s,X_{s}^{t,\eta},Y_{s}^{t,\eta},Z_{s}^{t,\eta},\sL_{(X _{s}^{t,\eta}, Y_{s}^{t,\eta},Z_{s}^{t,\eta})} \right)\dif s+Z_{s}^{t,\eta}\dif W_s, \hfill \\
Y_{T}^{t,\eta}=\Phi \left( X_{T}^{t,\eta},\sL_{X _{T}^{t,\eta}}\right), 
\end{array}
\label{equ12}	
\right .
\ee
and
\be \left \{ \begin{array}{l}
\dif Y_{s}^{t,x,\eta} \in A'\left( Y_{s}^{t,x,\eta} \right) \dif s-H\left( s,X_{s}^{t,x,\eta},Y_{s}^{t,x,\eta},Z_{s}^{t,x,\eta},\sL_{(X _{s}^{t,\eta}, Y_{s}^{t,\eta},Z_{s}^{t,\eta})}\right)\dif s+Z_{s}^{t,x,\eta}\dif W_s, \hfill \\
Y_{T}^{t,x,\eta}=\Phi \left( X_{T}^{t,x,\eta},\sL_{X _{T}^{t,\eta}}\right), 
\end{array}
\label{equ122}	
\right .
\ee 
where $A':\mR\mapsto 2^{\mR} $ is a maximal monotone operator and $H: \left[0,T \right] \times \mR^m\times \mR\times\mR^l\times\cM_2(\mR^m\times \mR\times\mR^l)\mapsto{\mR}, \Phi:\mR^m\times\cM_2(\mR^m)\mapsto{\overline{\cD(A')}}$ are Borel measurable. 

Assume:
\begin{enumerate}[(${\bf H}^1_{H,\Phi}$)] 
	\item The function $H$ is continuous in $(t,x,y,z,\vartheta)$, the function $\Phi$ is continuous in $(x,\mu)$ and they satisfy for $t\in[0,T], x\in\mR^m,y\in\mR,z\in\mR^l,\vartheta\in\cM_2(\mR^m\times \mR\times\mR^l), \mu\in\cM_2(\mR^m)$
	\ce
	{|{H(t,x,y,z,\vartheta)}|}+|\Phi(x,\mu)|\leq{L_{H,\Phi}(|{x}|+|{y}|+|{z}|+\|{\vartheta}\|_{2}+\|{\mu}\|_{2})},
	\de
	where $L_{H,\Phi}>0$ is a constant.
\end{enumerate}
\begin{enumerate}[(${\bf H}^2_{H,\Phi}$)] 
	\item The functions $H,\Phi$ satisfy for $t\in[0,T], x_1,x_2\in\mR^m, y_1,y_2\in\mR, z_1,z_2\in\mR^l,\vartheta_1,\vartheta_2\in\cM_2(\mR^m\times \mR\times\mR^l), \mu_1,\mu_2\in\cM_2(\mR^m)$,
	\ce
	&&\left|{H(t,x_1,y_1,z_1,\vartheta_1)}-{H(t,x_2,y_2,z_2,\vartheta_2)}\right| +\left| \Phi(x_1,\mu_1)-\Phi(x_2,\mu_2)\right| \\
	&\leq&L_{H,\Phi}(\left| {x_1-x_2}\right| +|{y_1-y_2}|+| {z_1-z_2}| +\rho(\vartheta_1,\vartheta_2)+\rho(\mu_1,\mu_2)). 
	\de
\end{enumerate}
\begin{enumerate}[(${\bf H}_{A'}$)]
\item $\cD(A')=\mR$ and $|A'^\circ(x)|\leq C(1+|x|), x\in\mR$, where $C>0$ is a constant.
\end{enumerate}

By Theorem \ref{EUni}, we know that Eq.(\ref{equ12}) has a unique solution denoted as $(Y_{\cdot}^{t,\eta},K_{\cdot}^{t,\eta}, Z_{\cdot}^{t,\eta})$. So, Eq.(\ref{equ122}) goes into a classical backward multivalued SDE. By Theorem \ref{bmsdeeu}, there exists a unique triple $(Y_{\cdot}^{t,x,\eta},K_{\cdot}^{t,x,\eta}, Z_{\cdot}^{t,x,\eta})$ satisfying Eq.(\ref{equ122}). About $Y_{\cdot}^{t,x,\eta}$ and $Z_{\cdot}^{t,x,\eta}$ we present the following estimate. 

\bp\label{ousoes}
Suppose that $({\bf H}^1_{H,\Phi})$-$({\bf H}^2_{H,\Phi})$ $({\bf H}_{A'})$ hold. Then there exists a constant $C>0$ such that 

(i) for any $t\in[0,T]$, $x\in\mR^m$, $\eta\in L^2(\tilde{\sF}_t,\mR^m)$, 
$$
\mE\left[\sup\limits_{s\in[t,T]}|Y_s^{t,x,\eta}|^2+\int_t^T|Z_s^{t,x,\eta}|^2\dif s|\tilde{\sF}_t\right]\leq C,
$$

(ii) for any $t\in[0,T]$, $x_1, x_2\in\mR^m$, $\eta_1, \eta_2\in L^2(\tilde{\sF}_t,\mR^m)$,
\ce
&&\mE\left[\sup\limits_{s\in[t,T]}|Y_s^{t,x_1,\eta_1}-Y_s^{t, x_2, \eta_2}|^2+\int_t^T|Z_s^{t,x_1,\eta_1}-Z_s^{t, x_2, \eta_2}|^2\dif s|\tilde{\sF}_t\right]\\
&\leq& C\left(|x_1-x_2|^2+\rho^2(\sL_{\eta_1},\sL_{\eta_2})\right).
\de
\ep
\begin{proof}
By $(ii)$ of Theorem \ref{diffesti}, we obtain $(i)$.

Next, we prove $(ii)$. The method is from \cite[Proposition 4.1]{L}. First of all, $( X^{t,x,\eta},Y^{t,x,\eta},\\ K^{t,x,\eta},Z^{t,x,\eta}) $ is independent of $\tilde{\mathscr{F}}_t$. Hence, it is independent of $\eta \in L^2\left( \tilde{\mathscr{F}}_t,\mathbb{R}^m \right)$. We consider $\left( X^{t,x,\eta},Y^{t,x,\eta},K^{t,x,\eta},Z^{t,x,\eta} \right) \left| _{x=\eta} \right. $. Combining the uniqueness of solutions for Eq.(\ref{equ12}) and Eq.(\ref{equ122}), with $X_{s}^{t,\eta}=X_{s}^{t,x,\eta}\left| _{x=\eta}=X_{s}^{t,\eta ,\eta} \right., s\in \left[ t,T \right]$, we get $\left( X^{t,\eta},Y^{t,\eta},K^{t,\eta},Z^{t,\eta} \right) =\left( X^{t,x,\eta},Y^{t,x,\eta},K^{t,x,\eta},Z^{t,x,\eta} \right) \left| _{x=\eta} \right. $. Besides, if $\overline{\eta }\in L^2\left( \tilde{\mathscr{F}}_t,\mathbb{R}^m \right)$, and the distribution of $\overline{\eta}$ is the same to that of $\eta$, then 
$$
\left( X^{t,\overline{\eta },\eta},Y^{t,\overline{\eta },\eta},K^{t,\overline{\eta },\eta},Z^{t,\overline{\eta },\eta} \right) :=\left( X^{t,x,\eta},Y^{t,x,\eta},K^{t,x,\eta},Z^{t,x,\eta} \right) \left| _{x=\overline{\eta }} \right. 
$$ 
and $\left( X^{t,\eta},Y^{t,\eta},K^{t,\eta},Z^{t,\eta} \right) $ have the same law.
So, given $\eta _{i},\eta _{i}^{'}\in L^2\left(\tilde{\mathscr{F}}_t,\mathbb{R}^m \right)$ with $\sL_{\eta _{i}}=\sL_{\eta^{'} _{i}}, i=1,2$, we consider the following equation,
\ce\left\{\begin{array}{l}
\dif Y_{s}^{t,\eta _{i}^{\prime},\eta _i}\in A'\left( Y_{s}^{t,\eta _{i}^{\prime},\eta _i} \right) \dif s-H\left( s,X_{s}^{t,\eta _{i}^{\prime},\eta _i},Y_{s}^{t,\eta _{i}^{\prime},\eta _i},Z_{s}^{t,\eta _{i}^{\prime},\eta _i},\mathscr{L}_{\left( X_{s}^{t,\eta _i},Y_{s}^{t,\eta _i},Z_{s}^{t,\eta _i} \right)} \right) \dif s+Z_{s}^{t,\eta _{i}^{\prime},\eta _i}\dif W_s\\
Y_{T}^{t,\eta _{i}^{\prime},\eta _i}=\varPhi \left( X_{T}^{t,\eta _{i}^{\prime},\eta _i},\mathscr{L}_{X_{T}^{t,\eta _i}} \right)	.
\end{array}
\right.
\de 
By Theorem \ref{diffesti} $(i)$ and $({\bf H}_{H,\varPhi}^{2})$, it holds that for any $\delta > 0$, there is a $\alpha> 0$ such that
\ce
&&\mathbb{E}\left[ \int_t^T{e^{\alpha \left( s-t \right)}}\left( \left| Y_{s}^{t,\eta _{1}^{\prime},\eta _1}-Y_{s}^{t,\eta _{2}^{\prime},\eta _2} \right|^2+\left\| Z_{s}^{t,\eta _{1}^{\prime},\eta _1}-Z_{s}^{t,\eta _{2}^{\prime},\eta _2} \right\| ^2 \right) \dif s \right] \no 
\\
&\leqslant& Ce^{\alpha T}\mathbb{E}\left[ \left| X_{T}^{t,\eta _{1}^{\prime},\eta _1}-X_{T}^{t,\eta _{2}^{\prime},\eta _2} \right|^2+\rho ^2\left( \mathscr{L}_{X_{T}^{t,\eta _1}},\mathscr{L}_{X_{T}^{t,\eta _2}} \right) \right] \no 
\\
&&+C\delta \mathbb{E}\Bigg[ \int_t^T e^{\alpha \left( s-t \right)}\Bigg( \left| X_{s}^{t,\eta _{1}^{\prime},\eta _1}-X_{s}^{t,\eta _{2}^{\prime},\eta _2} \right|^2\no\\
&&+\rho ^2\left( \mathscr{L}_{\left( X_{s}^{t,\eta _{1}^{\prime},\eta _1},Y_{s}^{t,\eta _{1}^{\prime},\eta _1},Z_{s}^{t,\eta _{1}^{\prime},\eta _1} \right)},\mathscr{L}_{\left( X_{s}^{t,\eta _{2}^{\prime},\eta _2},Y_{s}^{t,\eta _{2}^{\prime},\eta _2},Z_{s}^{t,\eta _{2}^{\prime},\eta _2} \right)} \right) \Bigg) \dif s\Bigg] \no \\
&\leqslant& C_{\alpha ,\delta}\left( \mathbb{E}\left[ \underset{s\in \left[ t,T \right]}{\sup}\left| X_{s}^{t,\eta _{1}^{\prime},\eta _1}-X_{s}^{t,\eta _{2}^{\prime},\eta _2} \right|^2+\rho ^2\left( \mathscr{L}_{X_{T}^{t,\eta _1}},\mathscr{L}_{X_{T}^{t,\eta _2}} \right) \right] \right) \no 
\\
&&+C\delta \int_t^T{e^{\alpha \left( s-t \right)}}\rho ^2\left( \mathscr{L}_{\left( X_{s}^{t,\eta _{1}^{\prime},\eta _1},Y_{s}^{t,\eta _{1}^{\prime},\eta _1},Z_{s}^{t,\eta _{1}^{\prime},\eta _1} \right)},\mathscr{L}_{\left( X_{s}^{t,\eta _{2}^{\prime},\eta _2},Y_{s}^{t,\eta _{2}^{\prime},\eta _2},Z_{s}^{t,\eta _{2}^{\prime},\eta _2} \right)} \right) \dif s \label{iequ2}.
\de 
Thus, by Lemma \ref{xest} and the definition of $\rho$, we get that
\ce 
&&\mathbb{E}\left[ \int_t^T{e^{\alpha \left( s-t \right)}}\left( \left| Y_{s}^{t,\eta _{1}^{\prime},\eta _1}-Y_{s}^{t,\eta _{2}^{\prime},\eta _2} \right|^2+\left\| Z_{s}^{t,\eta _{1}^{\prime},\eta _1}-Z_{s}^{t,\eta _{2}^{\prime},\eta _2} \right\| ^2 \right) \dif s \right] \no 
\\
&\leqslant& C_{\alpha ,\delta}\left( \mathbb{E}\left[ \left| \eta _{1}^{\prime}-\eta _{2}^{\prime} \right|^2+\rho ^2\left( \mathscr{L}_{\eta _1},\mathscr{L}_{\eta _2} \right) \right] \right) \no 
\\
&&+C\delta \int_t^T{e^{\alpha \left( s-t \right)}}\rho ^2\left( \mathscr{L}_{\left( X_{s}^{t,\eta _{1}^{\prime},\eta _1},Y_{s}^{t,\eta _{1}^{\prime},\eta _1},Z_{s}^{t,\eta _{1}^{\prime},\eta _1} \right)},\mathscr{L}_{\left( X_{s}^{t,\eta _{2}^{\prime},\eta _2},Y_{s}^{t,\eta _{2}^{\prime},\eta _2},Z_{s}^{t,\eta _{2}^{\prime},\eta _2} \right)} \right) \dif s\no
\\
&\leqslant& C_{\alpha ,\delta}\left( \mathbb{E}\left[ \left| \eta _{1}^{\prime}-\eta _{2}^{\prime} \right|^2+\rho ^2\left( \mathscr{L}_{\eta _1},\mathscr{L}_{\eta _2} \right) \right] \right) \no 
\\
&&+C\delta \mathbb{E}\left[ \int_t^T{e^{\alpha \left( s-t \right)}}\left( \left| Y_{s}^{t,\eta _{1}^{\prime},\eta _1}-Y_{s}^{t,\eta _{2}^{\prime},\eta _2} \right|^2+\left\| Z_{s}^{t,\eta _{1}^{\prime},\eta _1}-Z_{s}^{t,\eta _{2}^{\prime},\eta _2} \right\| ^2 \right) \dif s \right] .\label{iequ3}
\de 
Taking $\delta$ such that $C\delta < \frac{1}{2}$, we get 
\ce
&&\mathbb{E}\left[ \int_t^T{e^{\alpha \left( s-t \right)}}\left( \left| Y_{s}^{t,\eta _{1}^{\prime},\eta _1}-Y_{s}^{t,\eta _{2}^{\prime},\eta _2} \right|^2+\left\| Z_{s}^{t,\eta _{1}^{\prime},\eta _1}-Z_{s}^{t,\eta _{2}^{\prime},\eta _2} \right\| ^2 \right) \dif s \right] 
\\
&\leqslant& C_{\alpha ,\delta}\left( \mathbb{E}\left[ \left| \eta _{1}^{\prime}-\eta _{2}^{\prime} \right|^2+\rho ^2\left( \mathscr{L}_{\eta _1},\mathscr{L}_{\eta _2} \right) \right] \right).
\de

Besides, from the definition of $\rho$ and Lemma \ref{xest}, it follows that 
\ce
&&\int_t^T{\rho ^2\left( \mathscr{L}_{\left( X_{s}^{t,\eta _1},Y_{s}^{t,\eta _1},Z_{s}^{t,\eta _1} \right)},\mathscr{L}_{\left( X_{s}^{t,\eta _2},Y_{s}^{t,\eta _2},Z_{s}^{t,\eta _2} \right)} \right)}\dif s
\\
&=&\int_t^T{\rho ^2\left( \mathscr{L}_{\left( X_{s}^{t,\eta _{1}^{\prime},\eta _1},Y_{s}^{t,\eta _{1}^{\prime},\eta _1},Z_{s}^{t,\eta _{1}^{\prime},\eta _1} \right)},\mathscr{L}_{\left( X_{s}^{t,\eta _{2}^{\prime},\eta _2},Y_{s}^{t,\eta _{2}^{\prime},\eta _2},Z_{s}^{t,\eta _{2}^{\prime},\eta _2} \right)} \right)}\dif s
\\
&\leqslant& \mathbb{E}\left[\int_t^T \left( \left| X_{s}^{t,\eta _{1}^{\prime},\eta _1}-X_{s}^{t,\eta _{2}^{\prime},\eta _2} \right|^2+\left| Y_{s}^{t,\eta _{1}^{\prime},\eta _1}-Y_{s}^{t,\eta _{2}^{\prime},\eta _2} \right|^2+\left\| Z_{s}^{t,\eta _{1}^{\prime},\eta _1}-Z_{s}^{t,\eta _{2}^{\prime},\eta _2} \right\| ^2 \right) \dif s \right] 
\\
&\leqslant& C\left( \mathbb{E}\left[ \left| \eta _{1}^{\prime}-\eta _{2}^{\prime} \right|^2+\rho ^2\left( \mathscr{L}_{\eta _1},\mathscr{L}_{\eta _2} \right) \right] \right) .
\de 
So, taking the infimum on all $\eta _{1}^{\prime}, \eta _{2}^{\prime}$ with $\sL_{\eta _{i}^{\prime}}=\sL_{\eta _{i}}$ with $i=1,2$, we have 
$$
\int_t^T{\rho ^2\left( \mathscr{L}_{\left( X_{s}^{t,\eta _1},Y_{s}^{t,\eta _1},Z_{s}^{t,\eta _1} \right)},\mathscr{L}_{\left( X_{s}^{t,\eta _2},Y_{s}^{t,\eta _2},Z_{s}^{t,\eta _2} \right)} \right)}\dif s\leqslant C\rho ^2\left( \mathscr{L}_{\eta _1},\mathscr{L}_{\eta _2} \right) .
$$

Next, applying Theorem \ref{diffesti} $(ii)$ to Eq.(\ref{equ122}) with 
$$
g_i\left( s,y,z \right) :=H\left( s,X_{s}^{t,x,\eta _i},y,z,\mathscr{L}_{\left( X_{s}^{t,\eta _i},Y_{s}^{t,\eta _i},Z_{s}^{t,\eta _i} \right)} \right),  \quad \xi _i:=\varPhi \left( X_{T}^{t,x_i,\eta _i},\mathscr{L}_{X_{T}^{t,\eta _i}} \right), \quad i=1,2,
$$ 
by Lemma \ref{xest} we obtain that for all $\eta _1,\eta _2\in L^2\left( \tilde{\mathscr{F}}_t,\mathbb{R}^m \right) $, $x_1,x_2\in \mathbb{R}^m$, 
\ce
&&\mathbb{E}\left[ \underset{s\in \left[ t,T \right]}{\sup}\left| Y_{s}^{t,x_1,\eta _1}-Y_{s}^{t,x_2,\eta _2} \right|^2+\int_t^T{\left\| Z_{s}^{t,x_1,\eta _1}-Z_{s}^{t,x_2,\eta _2} \right\| ^2}\dif s\left| \tilde{\mathscr{F}}_t \right. \right] 
\\
&\leqslant& C\mathbb{E}\left[ \left| X_{T}^{t,x_1,\eta _1}-X_{T}^{t,x_2,\eta _2} \right|^2+\rho ^2\left( \mathscr{L}_{X_{T}^{t,\eta _1}},\mathscr{L}_{X_{T}^{t,\eta _2}} \right)+\int_t^T\left| X_{s}^{t,x_1,\eta _1}-X_{s}^{t,x_2,\eta _2} \right|^2\dif s \left| \tilde{\mathscr{F}} _t\right. \right] 
\\
&&+C\mathbb{E}\left[ \int_t^T{\rho ^2\left( \mathscr{L}_{\left( X_{s}^{t,\eta _1},Y_{s}^{t,\eta _1},Z_{s}^{t,\eta _1} \right)},\mathscr{L}_{\left( X_{s}^{t,\eta _2},Y_{s}^{t,\eta _2},Z_{s}^{t,\eta _2} \right)} \right) \dif s\left| \tilde{\mathscr{F}}_t \right.} \right] 
\\
&\leqslant& C\left( \left| x_1-x_2 \right|^2+\rho ^2\left( \mathscr{L}_{\eta _1},\mathscr{L}_{\eta _2} \right) \right).
\de
 The proof is complete.
\end{proof}

Based on the above proposition, we know that $Y_{\cdot}^{t,x,\eta}$ and $Z_{\cdot}^{t,x,\eta}$ depend on $\eta$ only through its distribution. Therefore, to strengthen the impression set
$$
Y_{s}^{t,x,\sL_{\eta}}:=Y_{s}^{t,x,\eta}, \quad Z_{s}^{t,x,\sL_{\eta}}:=Z_{s}^{t,x,\eta}, \quad s\in[t,T],
$$
and it holds that

(i) $(Y_{r}^{s,X_s^{t,x,\eta},\sL_{X_s^{t,\eta}}},Y_r^{s,X_s^{t,\eta}})=(Y_{r}^{t,x,\sL_{\eta}}, Y_{r}^{t,\eta}), \quad r\in[s,T], \mP.a.s.,$

(ii) $(Z_{r}^{s,X_s^{t,x,\eta},\sL_{X_s^{t,\eta}}},Z_r^{s,X_s^{t,\eta}})=(Z_{r}^{t,x,\sL_{\eta}}, Z_{r}^{t,\eta}), \quad r\in[s,T], \mP.a.s..$

Now set 
\be
u\left( t,x,\sL_{\eta} \right) :=Y_{t}^{t,x,\sL_{\eta}}, \quad \left( t,x,\sL_{\eta} \right) \in [0,T]\times\mR^m\times\cM_2(\mR^m),\label{equ13}
\ee 
and $u\left( t,x,\sL_{\eta} \right)$ is deterministic and satisfies
\ce
u(s,X_s^{t,x,\eta},\sL_{X_s^{t,\eta}})=Y_{s}^{s,X_s^{t,x,\eta},\sL_{X_s^{t,\eta}}}=Y_{s}^{t,x,\sL_{\eta}}.
\de

\medspace

Next, we introduce the following parabolic variation inequality (PVI for short) on $[0,T]\times\mR^m\times \mathcal{M}_2\left( \mathbb{R}^m \right)$:
\be \left \{ \begin{array}{l}
	\frac{\partial u\left( t,x,\sL_{\eta} \right)}{\partial t}+\mathcal{L}u\left( t,x,\sL_{\eta} \right) + H\left( t,x,u\left( t,x,\sL_{\eta} \right) ,\left( \nabla u\sigma \right) \left( t,x,\sL_{\eta} \right) ,\sL_{(\eta, u( t,\eta,\sL_{\eta}),\left( \nabla u\sigma \right) ( t,\eta,\sL_{\eta}))} \right)\\
	\in A'\left( u\left( t,x,\sL_{\eta} \right) \right), \\
	u\left( T,x,\sL_{\eta} \right) =\Phi \left( x,\sL_{\eta} \right) ,  \left( x,\sL_{\eta} \right) \in \mR^m\times \mathcal{M}_2\left( \mathbb{R}^m \right),
\end{array}
\label{equ14}	
\right.
\ee 
where 
\ce
\mathcal{L}u\left( t,x,\sL_{\eta} \right)&=&\left( b^i\partial _{x_i}u \right) \left( t,x,\sL_{\eta} \right) +\frac{1}{2}\left( \left( \sigma \sigma ^{\ast} \right) ^{ij}\partial _{x_ix_j}^{2}u \right) \left( t,x,\sL_{\eta} \right) \\
&&+\int_{\mR^m}\left( \partial _{\mu}u \right) _i\left( t,x,\sL_{\eta} \right)(y) b^i\left(t, y,\sL_{\eta} \right)\sL_{\eta}(\dif y) \\
&&+ \frac{1}{2}\int_{\mR^m}\partial _{y_i}\left( \partial _{\mu}u \right) _j\left( t,x,\sL_{\eta} \right)(y) \left( \sigma \sigma ^{\ast} \right) ^{ij}\left(t, y,\sL_{\eta} \right)\sL_{\eta}(\dif y) .
\de
Then we define viscosity solutions for PVI.(\ref{equ14}). For this, we introduce the following notations (cf.\cite{ADRI}):
\ce
A'_-\left( x\right):=\underset{x^{\prime}\rightarrow x, x^{\ast}\in A'\left( x^{\prime} \right)}{\lim\mathrm{inf}}x^{\ast}, \quad
A'_+\left( x\right):=\underset{x^{\prime}\rightarrow x, x^{\ast}\in A'\left( x^{\prime} \right)}{\lim\mathrm{sup}}x^{\ast}.
\de 

\bd\label{viscosity}
\begin{enumerate}[(i)]
\item we say that $u\in C([0,T]\times\mR^m\times\cM_2(\mR^m))$ is a viscosity subsolution of PVI.(\ref{equ14}) if $u(T,\cdot, \cdot)=\Phi(\cdot,\cdot)$ on $\mR^m\times\cM_2(\mR^m)$, and 
\ce
&&\frac{\partial \Psi \left( t,x,\sL_{\eta} \right)}{\partial t}+\mathcal{L}\Psi \left( t,x,\sL_{\eta} \right) \\
&&+ H\left( t,x,u\left( t,x,\sL_{\eta} \right) ,\left( \nabla \Psi \sigma \right) \left( t,x,\sL_{\eta} \right) ,\sL_{(\eta,u( t,\eta,\sL_{\eta}), ( \nabla \Psi \sigma)( t,\eta,\sL_{\eta}))}\right) \\
&\geq& A'_-\(u\left( t,x,\sL_{\eta}\right)\), 
\de 
whenever $\Psi\in C_b^{1,2,2}([0,T]\times\mR^m\times\cM_2(\mR^m))$ and 
$(t,x,\sL_{\eta})\in [0,T]\times\mR^m\times\cM_2(\mR^m) $ is a local maximum point of $u-\Psi$;
\item we say that $u\in C([0,T]\times\mR^m\times\cM_2(\mR^m))$ is a viscosity supersolution of PVI.(\ref{equ14}) if $u(T,\cdot, \cdot)=\Phi(\cdot,\cdot)$ on $\mR^m\times\cM_2(\mR^m)$, and 
\ce
&&\frac{\partial \Psi \left( t,x,\sL_{\eta} \right)}{\partial t}+\mathcal{L}\Psi \left( t,x,\sL_{\eta} \right)\\
&& + H\left( t,x,u\left( t,x,\sL_{\eta} \right) ,\left( \nabla \Psi \sigma \right) \left( t,x,\sL_{\eta} \right) ,\sL_{(\eta,u( t,\eta,\sL_{\eta}), ( \nabla \Psi \sigma)( t,\eta,\sL_{\eta}))}\right) \\
&\leq& A'_+\(u\left( t,x,\sL_{\eta} \right)\), 
\de 
whenever $\Psi\in C_b^{1,2,2}([0,T]\times\mR^m\times\cM_2(\mR^m))$ and 
$(t,x,\sL_{\eta})\in [0,T]\times\mR^m\times\cM_2(\mR^m)$ is a local minimum point of $u-\Psi$;
\item u is a viscosity solution of PVI.(\ref{equ14}) if it is both a viscosity subsolution and a viscosity supersolution of PVI.(\ref{equ14}).
\end{enumerate}
\ed

Now, we state the main result in this section.
\bt\label{Vico}
Assume that $({\bf H}^1_{b,\sigma})$-$({\bf H}^2_{b,\sigma})$, $({\bf H}^1_{H,\Phi})$-$({\bf H}^2_{H,\Phi})$ $({\bf H}_{A'})$ hold. Then the function $u(t,x,\sL_{\eta})$ defined by (\ref{equ13}) is a viscosity solution of PVI.(\ref{equ14}).
\et

To prove the above theorem, we make some preparation. First of all, for any $\varepsilon>0$, consider the following penalized equations 
\be
\left\{ \begin{array}{l}
	\mathrm{d}Y_{s,\varepsilon}^{t,\eta}=\[A'_{\varepsilon}\left( Y_{s,\varepsilon}^{t,\eta} \right)-H\left( s,X_{s}^{t,\eta},Y_{s,\varepsilon}^{t,\eta}, Z_{s,\varepsilon}^{t,\eta},\sL_{(X _{s}^{t,\eta},Y _{s,\varepsilon}^{t,\eta},Z_{s,\varepsilon}^{t,\eta})} \right)\]  \mathrm{d}s+Z_{s,\varepsilon}^{t,\eta}\mathrm{d}W_s,\\
	Y_{T,\varepsilon}^{t,\eta}=\Phi \left( X_{T}^{t,\eta},\sL_{X _{T}^{t,\eta}} \right), 
\end{array}
\label{equ15} 
\right. 
\ee
and
\be
\left\{ \begin{array}{l}
	\mathrm{d}Y_{s,\varepsilon}^{t,x,\eta}=\[A'_{\varepsilon}\left( Y_{s,\varepsilon}^{t,x,\eta} \right)-H\left( s,X_{s}^{t,x,\eta},Y_{s,\varepsilon}^{t,x,\eta},Z_{s,\varepsilon}^{t,x,\eta},\sL_{(X _{s}^{t,\eta},Y _{s,\varepsilon}^{t,\eta},Z_{s,\varepsilon}^{t,\eta})} \right)\]  \mathrm{d}s+Z_{s,\varepsilon}^{t,x,\eta}\mathrm{d}W_s,\\
	Y_{T,\varepsilon}^{t,x,\eta}=\Phi \left( X_{T}^{t,x,\eta},\sL_{X _{T}^{t,\eta}} \right), 
\end{array}
\label{equ155} 
\right. 
\ee 
where $A'_{\varepsilon}$ is the Yosida approximation of $A'$. So, by Theorem \ref{bmvsdeeu}, we know that under (${\bf H}^1_{H,\Phi}$)-(${\bf H}^2_{H,\Phi}$) Eq.(\ref{equ15}) has a unique solution denoted as $(Y_{\cdot,\varepsilon}^{t,\eta}, Z_{\cdot,\varepsilon}^{t,\eta})$. 
And Eq.(\ref{equ155}) is a classical backward SDE and has a unique solution denoted as $(Y_{\cdot,\varepsilon}^{t,x,\eta}, Z_{\cdot,\varepsilon}^{t,x,\eta})$ (cf. \cite[Theorem 4.1]{pp}). Besides, set $u_\e(t,x,\sL_\eta):=Y_{t,\varepsilon}^{t,x,\sL_\eta}:=Y_{t,\varepsilon}^{t,x,\eta}$, and by the similar deduction to that in the proof of Theorem \ref{EUni}, we obtain that $u_\e\rightarrow u$ as $\e\rightarrow 0$. Moreover, by \cite[Proposition 9.1 and Theorem 9.2]{L}, it holds that $u_\e(t,x,\sL_\eta)$ is continuous with respect to $(t,x,\sL_\eta)\in [0,T]\times\mR^m\times\cM_2(\mR^m)$ and a unique viscosity solution of the following parabolic partial differential equation:
\be \left \{ \begin{array}{l}
	\frac{\partial u_\e\left( t,x,\sL_{\eta} \right)}{\partial t}+\mathcal{L}u_\e\left( t,x,\sL_{\eta} \right) \\
	\quad + H\left( t,x,u_\e\left( t,x,\sL_{\eta} \right) ,\left( \nabla u_\e\sigma \right) \left( t,x,\sL_{\eta} \right) ,\sL_{(\eta, u_\e( t,\eta,\sL_{\eta}),\left( \nabla u_\e\sigma \right) ( t,\eta,\sL_{\eta}))} \right)\\
	=A'_\e\left( u_\e\left( t,x,\sL_{\eta} \right) \right), \\
	u_\e\left( T,x,\sL_{\eta} \right) =\Phi \left( x,\sL_{\eta} \right) ,  \left( x,\sL_{\eta} \right) \in \mR^m\times \mathcal{M}_2\left( \mathbb{R}^m \right).
\end{array}
\label{equ144}	
\right.
\ee

{\bf Proof of Theorem \ref{Vico}.} 

First of all, we prove that $u$ is a viscosity subsolution of PVI.(\ref{equ14}).

Let us take a $\Psi\in C_b^{1,2,2}([0,T]\times\mR^m\times\cM_2(\mR^m))$ such that $u-\Psi$ attains a local maximum in $(t, x,\sL_{\eta})\in [0,T]\times\mR^m\times\cM_2(\mR^m)$. Thus, there exists $\left( t_{\varepsilon},x_{\varepsilon},\eta _{\varepsilon}\right)$ such that, at least along a subsequence,

(i) $\left( t_{\varepsilon},x_{\varepsilon},\sL_{\eta _{\varepsilon}}\right)\rightarrow (t, x,\sL_\eta)$ as $\e\rightarrow 0$;

(ii) $u_\e-\Psi\leq u_\e( t_{\varepsilon},x_{\varepsilon},\sL_{\eta _{\varepsilon}})-\Psi( t_{\varepsilon},x_{\varepsilon},\sL_{\eta _{\varepsilon}})$ in a neighborhood of $( t_{\varepsilon},x_{\varepsilon},\sL_{\eta _{\varepsilon}})$ for any $\e>0$;

(iii) $u_\e( t_{\varepsilon},x_{\varepsilon},\sL_{\eta _{\varepsilon}})\rightarrow u(t, x,\sL_{\eta})$ as $\e\rightarrow 0$.

Since $u_\e$ is a viscosity subsolution of Eq.(\ref{equ144}), it holds that
\be
&&\frac{\partial \Psi \left( t_{\varepsilon},x_{\varepsilon},\sL_{\eta _{\varepsilon}} \right)}{\partial t}+\mathcal{L}\Psi \left( t_{\varepsilon},x_{\varepsilon},\sL_{\eta _{\varepsilon}} \right) + H\left( t_{\varepsilon},x_{\varepsilon},u_{\varepsilon}\left( t_{\varepsilon},x_{\varepsilon},\sL_{\eta _{\varepsilon}} \right) ,\left( \nabla \Psi \sigma \right) \left( t_{\varepsilon},x_{\varepsilon},\sL_{\eta _{\varepsilon}} \right),\vartheta_\e\right) 
\no\\
&\geqslant&A'_{\varepsilon}\left( u_{\varepsilon}\left( t_{\varepsilon},x_{\varepsilon},\sL_{\eta _{\varepsilon}} \right) \right).
\label{ineq}
\ee 
where $\vartheta_\e:=\sL_{(\eta _{\varepsilon},u_{\varepsilon}(t_{\varepsilon},\eta_{\varepsilon},\sL_{\eta _{\varepsilon}}),\left( \nabla \Psi \sigma \right) \left( t_{\varepsilon},\eta_{\varepsilon},\sL_{\eta _{\varepsilon}} \right))}$. Note that  
\ce
&&J_{\varepsilon}(u_{\varepsilon}\left( t_{\varepsilon},x_{\varepsilon},\sL_{\eta _{\varepsilon}} \right))\to u(t, x,\sL_{\eta}), \quad \e\rightarrow 0,\\
&& A'_{\varepsilon}\left( u_{\varepsilon}\left( t_{\varepsilon},x_{\varepsilon},\sL_{\eta _{\varepsilon}} \right) \right)\in A'\left( J_{\varepsilon}\left( u_{\varepsilon}\left( t_{\varepsilon},x_{\varepsilon},\sL_{\eta _{\varepsilon}} \right)\right) \right). 
\de
Thus, taking the lower limit on two sides of (\ref{ineq}) as $\e\rightarrow 0$, we obtain that
\ce
&&\frac{\partial \Psi \left( t,x,\sL_{\eta} \right)}{\partial t}+\mathcal{L}\Psi \left( t,x,\sL_{\eta} \right) + H\left( t,x,u\left( t,x,\sL_{\eta} \right) ,\left( \nabla \Psi \sigma \right) \left( t,x,\sL_{\eta} \right),\vartheta\right) \\
&\geqslant& A'_-(u\left( t,x,\sL_{\eta} \right)),
\de
where $\vartheta:=\sL_{(\eta,u( t,\eta,\sL_{\eta}), ( \nabla \Psi \sigma)( t,\eta,\sL_{\eta}))}$. So, by the definition, we know that $u$ is a viscosity subsolution of PVI.(\ref{equ14}).

By the same deduction to that for viscosity subsolutions of PVI.(\ref{equ14}), one can show that $u$ is a viscosity supersolution of PVI.(\ref{equ14}). The proof is complete.

\end{document}